\documentclass{amsart}

\usepackage{graphicx}
\usepackage{amssymb, amsmath,amsfonts, amsthm, hyperref, esint, mathtools, parskip}
\usepackage[english]{babel}
\usepackage[autostyle, english=american]{csquotes}
\usepackage[all,cmtip]{xy}
\MakeOuterQuote{"}

\newcommand{\norm}[1]{\left\lVert#1\right\rVert}

\theoremstyle{definition}
\newtheorem{theorem}{Theorem}[section]
\newtheorem{lemma}[theorem]{Lemma}
\newtheorem{proposition}[theorem]{Proposition}

\newtheorem{definition}[theorem]{Definition}
\newtheorem{example}[theorem]{Example}
\newtheorem{corollary}[theorem]{Corollary}

\title{Explicit Morphisms in the Galois-Tukey Category}
\author{David Philips}
\date{\today}

\begin{document}

\maketitle

\section{Introduction}

In 1874, Georg Cantor demonstrated that \(\mathfrak{c}\), the cardinality of the real numbers, is strictly greater than that of any countably infinite set \cite{cantor}. In other words, Cantor established that \(\aleph_0 < \mathfrak{c}\). This naturally leads to the question of whether there exists a cardinality between \(\aleph_0\) and \(\mathfrak{c}\). If no such intermediate cardinal exists, then \(\mathfrak{c} = \omega_1\), the smallest uncountable cardinal. This is precisely the statement of the Continuum Hypothesis (CH). Extending this idea, the Generalized Continuum Hypothesis (GCH) posits that for each ordinal number \(\alpha\), the cardinality of the power set of \(\omega_\alpha\) is equal to \(\omega_{\alpha+1}\).

In 1938, Kurt Gödel showed that CH is consistent with Zermelo-Fraenkel Set Theory combined with the Axiom of Choice (ZFC) by constructing a model of ZFC in which CH holds \cite{godel}. Building on Gödel's work, in 1963 Paul Cohen introduced the method of ``forcing'' which enabled him to construct a model of ZFC in which CH is false \cite{Cohen}. Together, Gödel’s and Cohen’s results imply that CH can neither be proven nor disproven using only the axioms of ZFC.

In this paper, we operate under the assumption that CH does \emph{not} hold. This allows for the existence of cardinalities between \(\omega_1\) and \(\mathfrak{c}\). These intermediate cardinalities are known as the "cardinal characteristics of the continuum." Each cardinal characteristic corresponds to a threshold at which some "nice" property of countable sets ceases to hold. For instance, the Baire Category Theorem ensures that a countable union of nowhere-dense sets cannot cover \(\mathbb{R}\) and this statement holds for $\aleph_0$ many nowhere dense sets. Therefore, we can ask how many more nowhere dense sets we need to add until the statement fails to generally hold. It turns out that the cardinality past which the Baire Category property fails is associated with the cardinal characteristic "$\mathrm{cov}(\mathcal{B})$." Similarly, \(\mathrm{add}(\mathcal{L})\) is the cardinal characteristic related to the property that the union of countably many sets of Lebesgue measure zero has measure zero. Later on, we will prove that $\mathrm{add}(\mathcal{L})\leq\mathrm{cov}(\mathcal{B})$.

Upon examining proofs related to cardinal characteristics, it becomes apparent that a substantial portion of them are very similar to each other. This redundancy inspired the construction of the "Galois-Tukey Category," first coined as $\mathbb{GT}$ by Peter Vojtáš \cite{Vojtáš}. Within $\mathbb{GT}$ analogous proofs are consolidated and their corresponding cardinal characteristics are treated as duals of one another. Furthermore, each inequality is associated with a morphism within the category. Thus, by leveraging \(\mathbb{GT}\) and the machinery of category theory, we can substantially simplify our study of cardinal characteristics while providing a cohesive framework for their analysis.

Although \(\mathbb{GT}\) is generally effective, it has certain limitations. Specifically, to work with a particular cardinal characteristic one must first define an appropriate "relation" for it. Although there is at least one trivial relation for any cardinal characteristic, some characteristics resist a working relation. Moreover, proving an inequality within \(\mathbb{GT}\) necessitates constructing a morphism between two relations. But, morphisms here require very stringent conditions to be met. There are even instances, in which a direct proof exists that one characteristic is of less than or equal cardinality than another, but constructing the corresponding morphism within \(\mathbb{GT}\) proves to be extremely difficult. Consequently, it is often unclear how to account for certain cardinal characteristics within $\mathbb{GT}$. 

My objective is to adapt some of the existing proofs, as well as give new proofs, of several established inequalities among cardinal characteristics so that their connection with \(\mathbb{GT}\) is made explicit. Through this process, I hope to highlight the strengths and limitations of the category. Specifically, by investigating which inequalities can or cannot be addressed within \(\mathbb{GT}\) and by comparing the complexity of proofs in \(\mathbb{GT}\) to direct proofs, we will get an idea of the extent to which the Galois-Tukey category serves as an effective framework for studying the cardinal characteristics of the continuum.

\section{The Galois-Tukey Category}
\subsection{Morphisms and Relations}\hphantom{.}

\begin{definition}\label{relation}
A triple $\mathbf{A}=(A_-,A_+,A)$ consisting of a set $A_-$, of "problems", another set $A_+$ of "solutions", and a binary relation $A\subseteq A_-\times A_+$, is called a \emph{relation}.  Here $xAy$ can be thought of as saying that "$y$ solves $x$."
\end{definition}

\begin{definition}\label{dual relation}
If $\mathbf{A}=(A_-,A_+,A)$ then the \emph{dual} of $\mathbf{A}$ is the relation $\mathbf{A^{\perp}} =(A_+, A_-,\neg {A^*})$ , where "$A^*$" is the converse of $A$. Here $(x,y)\in\neg{A^*}$ if and only if $(y,x)\not\in A$.
\end{definition}

Relations and their duals comprise the class of objects in \(\mathbb{GT}\). Since we will be utilizing relations to prove results about cardinality, we should have some notion that ties the two concepts together. This is done through the following definition. 

\begin{definition}\label{norm}
The \emph{norm} "$\norm{\mathbf{A}}$", of a relation $\mathbf{A}$, is the least cardinality of any subset $Y\subseteq A_+$, such that for each problem $x$ in $A_-$ there is a solution $y$ in $Y$ such that $xAy$. 
\end{definition}

\begin{example}\label{domination example}
For two functions $f,g\in\prescript{\omega}{\hphantom{.}}{\omega}$ we say $g<^* f$ (or $f$ \emph{dominates} $g$) if for all but finitely many $n\in\omega$, $g(n)<f(n)$. Equivalently, we say that there exists a point $n_0\in\omega$ such that for every $n>n_0$, $g(n)<f(n)$. Although domination is transitive, it does not define a total ordering. For example, if $f\in\prescript{\omega}{\hphantom{.}}{2}$ is the function that only returns $1$ and $g\in\prescript{\omega}{\hphantom{.}}{2}$ is a function that alternates between $0$ and $1$, neither dominates the other. 

Let $\mathfrak{D}$ be the relation $(\prescript{\omega}{\hphantom{.}}{\omega},\prescript{\omega}{\hphantom{.}}{\omega},<^*)$, then we define $\norm{\mathfrak{D}}:=\mathfrak{d}$. Likewise, $\norm{\mathfrak{D^\perp}}=\norm{(\prescript{\omega}{\hphantom{.}}{\omega},\prescript{\omega}{\hphantom{.}}{\omega},\not>^*)}:=\mathfrak{b}$. We call $\mathfrak{d}$ the "dominating number" and $\mathfrak{b}$ the "bounding number."

\begin{example}\label{alternate domination example}

An \emph{interval partition} is a partition of \(\omega\) into infinitely many finite intervals. We say an interval partition $I$, \emph{dominates} another interval partition  $J$, if for all but finitely many  $j_k\in J$ there exists some $i_n\in I$ such that $j_k\subseteq i_n$. 
We denote "\(\mathrm{IP}\)" as the set of all interval partitions of $\omega$.

Alternatively, if we define $\mathfrak{D'}:=(\mathrm{IP},\mathrm{IP},\text{dominated by})$, it can easily be shown that $\mathfrak{d}=\norm{\mathfrak{D'}}$ and $\mathfrak{b}=\norm{\mathfrak{D'^{\perp}}}$.
\end{example}

These two examples illustrate that a single cardinal can be associated with multiple relations. This flexibility is advantageous when proving inequalities between cardinal characteristics because we can always pick which relation according to what suits our needs.
\end{example}

\begin{example}\label{splitting example}
For $x,y\in[\omega]^\omega,$ we say $x$ \emph{splits} $y$ if both $x\cap y$ and $y\setminus x$ are infinite. 

Let $\mathfrak{R}$ be the relation $(\mathcal{P}(\omega),[\omega]^\omega,\text{does not split})$, then $\norm{\mathfrak{R}}:=\mathfrak{r}.$ Likewise, $\norm{\mathfrak{R^\perp}}=\norm{([\omega]^\omega,\mathcal{P}(\omega),\text{split by})}:=\mathfrak{s}$. We call $\mathfrak{r}$ the "reaping number" and $\mathfrak{s}$ the "splitting number."
\end{example}

\begin{example}\label{coloring example}
For an infinite set $Q\in[\omega]^\omega$ and a $2$-coloring $\pi:[\omega]^n\rightarrow 2$, we take $\pi H Q$ to mean that there is a set $K\in fin(\omega)$, such that for every $q\in[Q\setminus K]^n$, $\pi(q)$ is constant. In other words, we say that $Q$ is \emph{almost homogeneous} for $\pi$.

Letting $P_n$ represent the set of all $2$-colorings $\pi:[\omega]^n\rightarrow 2$, we define the relation  $\mathfrak{Hom_n}:= (P_n,[\omega]^\omega,H)$. Then $\norm{\mathfrak{Hom_n}}:=\mathfrak{hom_n}$ and $\norm{\mathfrak{Hom_n^\perp}}=\norm{([\omega]^\omega,P_n,\neg H^*)}:=\mathfrak{par_n}$. $\mathfrak{hom_n}$ is called the "homogeneity number" and $\mathfrak{par_n}$ is called the "partition number."
\end{example}

All categories are composed of objects and morphisms. Having defined the objects of $\mathbb{GT}$, we will now specify its morphisms.

\begin{definition}\label{morphism}
A \emph{morphism} between two relations, $\mathbf{A}=(A_-,A_+,A)$ and $\mathbf{B}=(B_-,B_+,B)$, is a pair of functions $\varphi=(\varphi_-:B_-\rightarrow A_-,\varphi_+:A_+\rightarrow B_+)$ such that for every $b\in B_-$ and $a\in A_+$, if $\varphi_- (b)Aa$ then $b B\varphi_+ (a)$. As a shorthand, we let $\varphi:\mathbf{A}\rightarrow\mathbf{B}$ denote a morphism $\varphi$ from $\mathbf{A}$ to $\mathbf{B}.$
\end{definition}

\begin{figure}[h]
    \centering
    \includegraphics[width=0.4\textwidth]{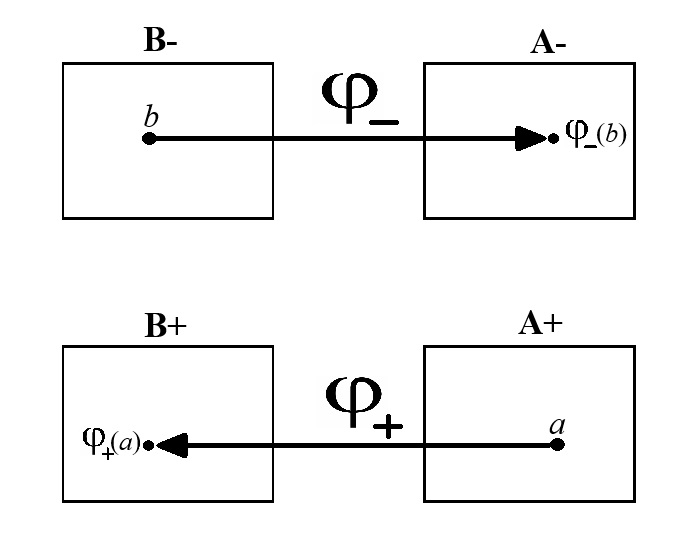}
    \caption{An Illustration of a Morphism}
\end{figure}

Notice that in defining a morphism $\varphi=(\varphi_-,\varphi_+)$ from $\mathbf{A}$ to $\mathbf{B}$, we automatically get a dual morphism $\varphi^\perp=(\varphi_+,\varphi_-)$ from $\mathbf{B^\perp}$ to $\mathbf{A^\perp}$. 

Defining a relation for a specific cardinal is already cumbersome. Moreover, it is not always clear how to define a relation whose dual corresponds to a different cardinal. So one might wonder why we choose to work within such a category. The reason lies in the following, simple but powerful, theorem.

\begin{theorem}\label{morphism theorem}
If there exists a morphism $\varphi:\mathbf{A}\rightarrow \mathbf{B}$ then $\norm{\mathbf{A}}\geq\norm{\mathbf{B}}$ and $\norm{\mathbf{A^{\perp}}}\leq\norm{\mathbf{B^{\perp}}}$.
\end{theorem}

\begin{proof}
The latter inequality follows immediately after applying the former inequality to the dual morphism $\varphi^{\perp}: \mathbf{B^{\perp}}\rightarrow \mathbf{A^{\perp}}$, so we will only prove the first. 

Suppose $\varphi:\mathbf{A}\rightarrow \mathbf{B}$ is a morphism from $\mathbf{A}=(A_-,A_+,A)$ to $\mathbf{B}=(B_-,B_+,B)$. Let $X\subseteq A_+$ be a cofinal set with cardinality $\norm{\mathbf{A}}$. Since $|\varphi_+(X)|\leq|X|$, it suffices to show that $\varphi_+(X)$ is cofinal in $B_+$. To this end,  let $b\in B_-$ be arbitrary and suppose $\varphi_-(b)Ay$, for some $y\in A_+$. By the definition of a morphism, $b B\varphi_+(y)$.
\end{proof}

\begin{definition}\label{sigma-relation}
Given a relation $\mathbf{A}=(A_-,A_+,A)$, we can define a corresponding sigma-relation, $$\mathbf{A_\sigma}=(\prescript{\omega}{\hphantom{.}}{A_-,A_+,A_\sigma}).$$ $xA_\sigma y$ means that $x(n)Ay$ for every $n\in\omega$. $\norm{\mathbf{A_\sigma}}$ represents the least cardinality of a solution set $Y\subseteq A_+$ such that any countable set of problems $X\subseteq A_-$ can be solved by at least one $y\in Y$.
\end{definition}

 It is always the case that $\norm{\mathbf{A_\sigma}}\geq\norm{\mathbf{A}}$, as a morphism can be trivially constructed by allowing $\varphi_-:A_-\rightarrow\prescript{\omega}{\hphantom{.}}{A_-}$ to be the function that associates elements $a\in A_-$ with any countable set including $a$ within it and letting $\varphi_+$ be the identity. For many cardinal characteristics, it remains an open question as to whether or not this inequality can be made strict. 

\begin{example}\label{sigma example}
Notice $\mathfrak{R}_\sigma=(^\omega \mathcal{P}(\omega), [\omega]^\omega, \text{does not split})$. Then 
$\|\mathfrak{R}_\sigma\|:= \mathfrak{r}_\sigma\geq\mathfrak{r}$ and $\|\mathfrak{R}_\sigma^\perp\| = \mathfrak{s}$. Whether or not $\mathfrak{r}_\sigma>\mathfrak{r}$ is consistent with ZFC is an open problem. 
\end{example}

From the above examples, we learn that norms of relations can represent certain cardinal characteristics and from the preceding theorem, we learn that morphisms between relations induce inequalities between these norms. Since our goal is to study inequalities among cardinal characteristics, it is unsurprising that much of this thesis will focus on explicitly constructing morphisms between various relations. Before moving on to more general categorical questions we will give one more useful result regarding morphisms.

\begin{theorem}\label{morphism theorem 2}
Suppose $\mathbf{A} = (A_-, A_+, A)$ and $\mathbf{B} = (B_-, B_+, B)$ are two relations and $\kappa\geq\aleph_0$ is a cardinal.
\begin{enumerate}
    \item $\|\mathbf{A}\| \leq \kappa$ if and only if there exists a morphism from $(\kappa, \kappa, =)$ to $\mathbf{A}$.
    \item If $\|\mathbf{A}\| = |A_+| = \kappa$, then there exists a morphism from $\mathbf{A}$ to $(\kappa, \kappa, <)$.
    \item If $\|\mathbf{A}^{\perp}\| = |A_-| = \kappa$, then there exists a morphism from $(\kappa, \kappa, <)$ to $\mathbf{A}$.
    \item If $\|\mathbf{A}\| = |A_+| = \norm{\mathbf{B}^\perp} = |B_-| \geq \aleph_0$, then there exists a morphism from $\mathbf{A}$ to $\mathbf{B}$.
\end{enumerate}
\pagebreak
\begin{proof}\phantom{.}
 \begin{enumerate}
     \item We first construct a morphism $\varphi:(\kappa,\kappa,=)\rightarrow\mathbf{A}$. Assuming that $\norm{\mathbf{A}}\leq \kappa$, we can let $\varphi_+:\kappa\rightarrow A_+$ surject onto a cofinal set of solutions within $A_+$. For an arbitrary $a\in A_-$ we will let $\varphi_-(a)$ be equal to some $\delta<\kappa$, such that $\varphi_+(\delta)$ is a solution for $a$. To verify this is a morphism, we must check that $\varphi_{-}(a)=\delta$ implies $aA\varphi_{+}(\delta)$. This is immediate by the construction of $\varphi_-$. The other direction is obvious because of Theorem \ref{morphism theorem}.
     \item Let $\varphi_+:A_+\rightarrow\kappa$ be an injection. For any $\delta<\kappa$, we can define the set $X_\delta:=\{a\in A_+:\varphi_+(a)\leq\delta\}$. Since $|X_\delta|<\kappa$, there exists at least one $a_\delta\in A_-$ which is not solved by any element of $X_\delta$. We will let $\varphi_{-}$ be the function associating each $\delta$ with some corresponding $a_\delta$. Assuming that $a_\delta Aa$, for some $a\in A_-$, we get $\varphi_{+}(a)>\delta$.
     \item If we apply (2) to the dual relation $\mathbf{A}^\perp$, we get a morphism from $\mathbf{A}^\perp$ into $(\kappa,\kappa,<)$
     \item From (2) we get a morphism $\varphi_1:\mathbf{A}\rightarrow(\kappa,\kappa,<)$ and from (3) we get a morphism $\varphi_2:(\kappa,\kappa,<)\rightarrow\mathbf{B}$. Then $\varphi:=\varphi_1\circ\varphi_2$, is a morphism from $\mathbf{A}$ to $\mathbf{B}$. The fact that we can compose morphisms is proved in Proposition \ref{category theorem}.
 \end{enumerate}   
\end{proof}
\end{theorem}
\vspace{-0.3cm}

\subsection{The Category $\mathbb{GT}$}\hphantom{.}

\begin{proposition}\label{category theorem}
$\mathbb{GT}$ is a category.
\begin{proof}\phantom{.}
\begin{enumerate}
    \item \emph{Identity Morphisms:} 
    Let $\mathbf{A}\in\text{Obj}(\mathbb{GT})$ be arbitrary. We will define the identity morphism as $\mathrm{id}_{\mathbf{A}} := ( \mathrm{id}_{A_-}, \mathrm{id}_{A_+} ),$
    where $\mathrm{id}_{A_-}$ is the identity map on $A_-$, and $\mathrm{id}_{A_+}$ is the identity map on $A_+$. For any $x \in A_-$ and $y \in A_+$, if $\mathrm{id}_{A_-}(x) A y$, then $x A \:\mathrm{id}_{A_+}(y)$. Hence, $\mathrm{id}_{\mathbf{A}} : \mathbf{A} \to \mathbf{A}$ is a morphism.

  \item \emph{Composition of Morphisms:} Let $\varphi=(\varphi_-,\varphi_+):\mathbf{A}\to\mathbf{B}$ and $\psi=(\psi_-,\psi_+):\mathbf{B}\to\mathbf{C}$. We define $(\psi\circ\varphi)_-:=\varphi_-\circ\psi_-$ and $(\psi\circ\varphi)_+:=\psi_+\circ\varphi_+$. Suppose $(\psi\circ\varphi)_-(c)\,A\,a$, meaning $\varphi_-(\psi_-(c))\,A\,a$. Since $\varphi$ is a morphism, we get $\psi_-(c)\,B\,\varphi_+(a)$, and because $\psi$ is a morphism, $c\,C\,\psi_+(\varphi_+(a))$. Hence, $cC(\psi\circ\varphi)_+(a)$, and so $\psi\circ\varphi$ is a morphism.

\item \emph{Associativity of Composition:} This follows from the associativity of set functions. 
\end{enumerate}
\end{proof}
\end{proposition}

We can also ask whether or not $\mathbb{GT}$ contains zero objects. The answer is a resounding "no." In fact, $\mathbb{GT}$ contains neither initial nor final objects.

\begin{proposition}\label{terminal object theorem}
There are no zero objects in $\mathbb{GT}.$
\end{proposition}

\begin{proof}
    Suppose, for the sake of contradiction, that $\mathbf{A}:=(A_-,A_+,A)$ was initial. Let $X$ be a set of strictly greater cardinality than $A_+$ and define the relation $\mathbf{X}:=(X,X,=)$. Since $\mathbf{A}$ is initial, there exists a morphism $\varphi:\mathbf{A}\rightarrow \mathbf{X}$. By Theorem \ref{morphism theorem}, this means that $A_+\geq\norm{\mathbf{A}}\geq\norm{\mathbf{X}}=|X|>A_+$.

    If we instead assumed that $\mathbf{A}$ was final, this would imply there is a morphism $\varphi:\mathbf{B}\rightarrow \mathbf{A}$, for any arbitrary relation $\mathbf{B}$. By applying the above argument to $\mathbf{A^\perp}$ we end up with the same contradiction.
\end{proof}

Although $\mathbb{GT}$ lacks zero objects, there is still structure to be found. In particular, we can define finite products and co-products.

\begin{definition}\label{product/coproduct}
The \emph{product} $\mathbf{A\times B}$ of two relations is defined as the relation $(A_{-}\sqcup B_{-},A_{+}\times B_{+},C)$.  $(a,b) C(x,y)$ means that if $a\in A_{-}$ and $b=0$, then $aAx$. If instead $b\in B_{-}$ and $a=1$, then $(a,b) C(x,y)$ means $bBy$. The \emph{coproduct} $\mathbf{A+B}$ is defined as $(\mathbf{A^{\perp}}\times \mathbf{B}^{\perp})^{\perp}$.
\end{definition}
\begin{theorem}\label{relation operations}\phantom{.}
\begin{enumerate}
    \item  $ \| \mathbf{A} \times \mathbf{B} \|  =  \max \{ \| \mathbf{A} \|, \| \mathbf{B} \| \}. $ 
    \item $ \| \mathbf{A} + \mathbf{B} \|  =  \min \{ \| \mathbf{A} \|, \| \mathbf{B} \| \}.
$  
\end{enumerate}
\end{theorem}
\begin{proof}\phantom{.}
We will only prove (1) as (2) can be proven with a dual argument.  By definition, a cofinal set  $X\subseteq A_+\times B_+$ contains a solution for every problem $(a,0)\in A_-\sqcup B_-$. By associating each pair $(a,0)$ with its corresponding element $a\in A_-$, we see that $\|\mathbf{A} \times \mathbf{B}\|\geq\norm{\mathbf{A}}$. Likewise, by associating each pair $(1,b)\in A_-\sqcup B_-$ with its corresponding element $b\in B_-$, we get that $\|\mathbf{A} \times \mathbf{B}\|\geq\norm{\mathbf{B}}$.

To show that $\norm{\mathbf{A}\times\mathbf{B}}\leq \max{\{\norm{\mathbf{A}},\norm{\mathbf{B}}\}}$,  we separately consider the case where both $\norm{\mathbf{A}}$ and $\norm{\mathbf{B}}$ are finite and the case where at least one of the two is infinite. Starting with the infinite case, let $S_A\subseteq A_+$ be a cofinal set with $|S_A|=\norm{\mathbf{A}}$ and let $S_B\subseteq B_+$ be a cofinal set with $|S_B|=\norm{\mathbf{B}}$. Fix $p\in B_+$ and $q\in A_+$ and define: $$S:=\{(a,p):a\in S_A\}\cup\{(q,b):b\in S_B\}\subseteq A_+\times B_+.$$

By the fact that $S$ is a solution set for $A_-\sqcup B_-$ and properties of infinite cardinal arithmetic, $$\norm{\mathbf{A}\times\mathbf{B}}\leq|S|\leq|S_A|+|S_B|=\max\{\norm{\mathbf{A}},\norm{\mathbf{B}}\}.$$

If both $\norm{\mathbf{A}}$ and $\norm{\mathbf{B}}$ are finite, let $|S_A|=n$ and $|S_B|=m$. Without loss of generality, assume that $n\geq m$. Define : $$T := \{(x_i, y_i) : 1 \leq i \leq m\} \cup \{(x_j, y_m) : m+1 \leq j \leq n\}.
$$
If $(a,0)\in A_-\sqcup B_-$ is a problem, then there must be some $x_i\in S_A$ such that $aAx_i$. If $1\leq i\leq m$ then, $(x_i,y_i)$ is a solution. If $m<i\leq n$ then $(x_i,y_m)$ is a solution. If $(1,b)\in A_-\sqcup B_-$ is a problem, then  $(x_i,y_i)$ is a solution. This shows that $T$ is a solution set for $A_-\sqcup B_-$. Thus, $$\norm{\mathbf{A}\times\mathbf{B}} \leq |T| \leq \max\{\norm{\mathbf{A}},\norm{\mathbf{B}}\}.$$
\end{proof}
\pagebreak
Beyond the product and co-product, there are several other ways we can combine two relations. 

\begin{definition}\label{alternate definitions}{For two relations $\mathbf{A}=(A_-,A+,A)$ and $\mathbf{B}=(B_-,B+,B)$}:
\begin{enumerate}
    \item The \emph{conjunction} is the relation $\mathbf{A\land B}:=(A_-\times B_-,A_+\times B_+,K)$, where $(x,y)K(a,b)$ means $xAa$ and $yBb$
    \item The \emph{sequential composition} is the relation  $\mathbf{A;B}:=(A_-\times{}^{A_+}B_-,A_+\times B_+,S)$. Where $(x,f)S(a,b)$ means $xAa$ and $f(a)Bb$.
    \item  The \emph{dual sequential composition} is the relation $\mathbf{A;^*B:=(A^{\perp};B^{\perp}})^{\perp}$.
\end{enumerate}
\end{definition}

As shown below, properties analogous to Theorem \ref{relation operations} remain true for relations defined above. The proofs of which are similar to the argument given in Theorem \ref{relation operations} and are therefore omitted. 

\begin{theorem}\label{extra relation operations}
\end{theorem}
\begin{enumerate}
    \item$ \max \{ \| \mathbf{A} \|, \| \mathbf{B} \| \}  \leq  \| \mathbf{A} \land \mathbf{B} \|  \leq  \| \mathbf{A} \| \cdot \| \mathbf{B} \|. $  \\
    \item $ \| \mathbf{A};\mathbf{B}  \|  =  \| \mathbf{A} \| \cdot \| \mathbf{B} \|. $  \\
    \item$ \| \mathbf{A};^* \mathbf{B} \|  =  \min \{ \| \mathbf{A} \|, \| \mathbf{B} \| \}. $
\end{enumerate}

In the infinite case, maximums and products are the same. For (1) this means that $ \max \{ \| \mathbf{A} \|, \| \mathbf{B} \| \}  =  \| \mathbf{A} \land \mathbf{B} \|$ and for (2) this means that  $ \| \mathbf{A};\mathbf{B}  \|=\max \{ \| \mathbf{A} \|, \| \mathbf{B} \| \}$. 
\section{Cichoń's Diagram}

Cichoń's diagram is a crucial tool in the study of cardinal characteristics of the continuum. It elegantly organizes several important cardinal invariants that measure different properties of the real line. As depicted here, the diagram is backward in comparison to its traditional orientation. That is, all the below arrows are reversed. Usually, an arrow from one cardinal characteristic $X$ to another characteristic $Y$ represents the fact that $X\leq Y$. For our purposes,  we take an arrow from $X$ to $Y$ to mean that there exists a morphism $\varphi:\mathbf{X}\rightarrow\mathbf{Y}$. Where $\mathbf{X}$ and $\mathbf{Y}$ represent relations whose norms are $X$ and $Y$, respectively. Thus, the reason we reverse arrows boils down to Theorem \ref{morphism theorem}.
\begin{figure}[h]
    \centering
\[
\xymatrix@R-0.5pc@C-0.8pc{
&
\operatorname{cov}(\mathcal{L}) \ar[dd] &
\operatorname{non}(\mathcal{B}) \ar[d]\ar[l] &
\operatorname{cof}(\mathcal{B}) \ar[d]\ar[l] &
\operatorname{cof}(\mathcal{L}) \ar[l]\ar[dd] &
\ar[l]\mathfrak{c}
\\
&& \mathfrak{b} \ar[d] & \mathfrak{d} \ar[l]\ar[d]&
\\
\omega_1
&
\operatorname{add}(\mathcal{L})\ar[l]&
\operatorname{add}(\mathcal{B}) \ar[l] &
\operatorname{cov}(\mathcal{B}) \ar[l] &
\operatorname{non}(\mathcal{L}) \ar[l] &
}
\]
\caption{Cichoń's Diagram}
    \label{cichon}
\end{figure}
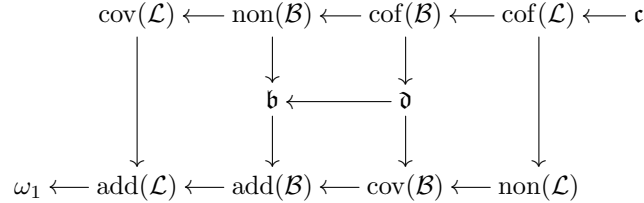

One of the primary reasons Cichoń's diagram is significant is because of its ability to unify seemingly separate areas within set theory. Cardinal characteristics related to Lebesgue measure, such as $\mathrm{cov}(\mathcal{L})$ and $\mathrm{add}(\mathcal{L})$, find their counterparts in the realm of Baire category through invariants like $\mathrm{cov}(\mathcal{B})$ and $\mathrm{add}(\mathcal{B})$. This unification reveals connections between measure-theoretic and topological properties of $\mathbb{R}$. Moreover, Cichoń's diagram plays a role in independence results. Different models of set theory, achieved through various forcing extensions or by adopting alternative axioms beyond ZFC, can result in distinct versions of the diagram in which the relative sizes of the cardinal characteristics differ \cite{switzer}. However, the version presented here represents the standard configuration within ZFC. This canonical diagram captures the foundational ZFC-provable inequalities and serves as a baseline for our analysis.

The main objective of this section will be to explicitly construct each morphism presented in Figure \ref{cichon}. To achieve this, we will start by formally defining each of the depicted cardinal characteristics. Following these definitions, we will introduce relations whose norms correspond to each defined cardinal. After all of this, we will construct the relevant morphisms. 

\subsection{Bounding and Dominating}\hphantom{.}

As defined in examples \ref{domination example} and \ref{alternate domination example}, the two cardinals from Cichoń's diagram that we start with are $\mathfrak{b}$ and $\mathfrak{d}$.

\begin{theorem}\label{b<d}
$\mathfrak{b}\leq\mathfrak{d}$.
\begin{proof}
By Example \ref{domination example}, we associate $\mathfrak{d}$ with the norm of the relation $(\prescript{\omega}{\hphantom{.}}{\omega},\prescript{\omega}{\hphantom{.}}{\omega},<^*)$ and $\mathfrak{b}$ with the norm of the dual relation $(\prescript{\omega}{\hphantom{.}}{\omega},\prescript{\omega}{\hphantom{.}}{\omega},\not>^*)$.

Let both $\varphi_-,\varphi_+:\prescript{\omega}{\hphantom{.}}{\omega}\rightarrow\prescript{\omega}{\hphantom{.}}{\omega}$ be the identity function. Then $\varphi_-(f)<^*g$ implies $f<^*\varphi_+(g).$ Since domination is antisymmetric, $f\not>^*\varphi_+(g)$.
\end{proof}
\end{theorem}

After having proved Theorem \ref{b<d}, we need to establish that $\mathfrak{b}$ is uncountable. Doing so within $\mathbb{GT}$ requires a morphism from $(\prescript{\omega}{\hphantom{.}}{\omega},\prescript{\omega}{\hphantom{.}}{\omega},\not>^*)$ into a relation whose norm is $\omega_1$. The most natural candidate for such a relation would be something like $(\omega_1,\omega_1,=)$. We now show that no such morphism exists. In fact no relation $\mathbf{A}=(A_-,A_+,A)$ with $|A_-|<\mathfrak{d}$ will work.

\begin{proposition}\label{noMorph}
Let $\mathbf{A}:=(A_-,A_+,A)$ be a relation such that  $\norm{\mathbf{A}}>1$. If $\varphi:(\prescript{\omega}{\hphantom{.}}{\omega},\prescript{\omega}{\hphantom{.}}{\omega},\not>^*)\rightarrow\mathbf{A}$ is a morphism, then $|A_-|\geq\mathfrak{d}$.
\end{proposition}
\begin{proof}
For contradiction, suppose $|A_-|<\mathfrak{d}.$ Define $\mathcal{D}:=\{\varphi_-(a):a\in A_-\}\subseteq \prescript{\omega}{\hphantom{.}}{\omega}$. Since $\mathcal{D}$ is not a dominating family, let $h\in \prescript{\omega}{\hphantom{.}}{\omega}$ be a function not dominated by any element of $\mathcal{D}$. By the definition of a morphism, for every $a\in A_-$, $aA\varphi_+(h)$. This implies $\norm{\mathbf{A}}=1$. By contradiction, $\mathcal{D}$ is a dominating family. But this is impossible since $|A_-|<\mathfrak{d}$. 
\end{proof}

 It should be noted that we can easily get the above result by considering the dual morphism $\varphi^\perp$, but the details of this proof will serve as an interesting comparison later on with Proposition \ref{noMorph2} and Theorem \ref{w<par}. Moreover, Proposition \ref{noMorph} can be generalized to similar relations with norm $\mathfrak{b}$. For example, the alternate construction given in Example \ref{alternate domination example} will also fail because there are morphisms to and from it and the standard relation for $\mathfrak{b}$. Although it may be possible to construct a relation and a corresponding morphism that work in this scenario, it is likely difficult. As such a relation could not have a morphism from the standard witnesses for $\mathfrak{b}$ into it and its problem set must have strictly greater cardinality than $\mathfrak{d}$. This is all to say that giving a direct proof is much easier.

\begin{proposition}\label{w<b}
$\omega_1\leq\mathfrak{b}$   
\end{proposition}
\begin{proof}
Let $\mathcal{G}:=\{g_n\in\prescript{\omega}{\hphantom{.}}{\omega}: n\in\omega\}$ be an arbitrary countable family of functions. For every $k\in\omega$ define the function $f\in\prescript{\omega}{\hphantom{.}}{\omega}$ as: $$f(k)=\bigcup\{g_i(k):i\in k\}.$$
Since $f$ dominates each $g_i\in\mathcal{G}$, we can say that at least one function dominates any countable family of functions.
\end{proof}
\subsection{Ideals}\hphantom{.}

\begin{definition}\label{IDEAL}
Given a set \( X \), an \emph{ideal} \( I \) on \( X \) is a nonempty subset of \( \mathcal{P}(X) \) such that:
\begin{enumerate}
    \item If \( A \in I \) and \( B \subseteq A \), then \( B \in I \).
    \item If \( A, B \in I \), then \( A \cup B \in I \).
    \item \( X \notin I \), ($I$ is proper).
\end{enumerate}
\end{definition}

It should be clear that ideals are dual to filters. Less obviously, recall that $\mathcal{P}(X)$ can be transformed into a commutative ring by defining addition with symmetric differences, multiplication by intersections, and letting the additive identity be the empty set. Under these conditions $\mathcal{P}(X)$ is a Boolean ring, meaning that each of its elements is idempotent under multiplication. In this context, subsets $I\subset\mathcal{P}(X)$ conforming to the conditions of Definition \ref{IDEAL} are ideals in the ring theoretic sense.
\begin{definition}\label{ideal definition}
    Let $\mathcal{I}$ be a proper ideal of subsets of a set $X$ which contains all of its singletons.

\begin{enumerate}
    \item The \emph{additivity} of $\mathcal{I}$, $\mathrm{add}(\mathcal{I})$, is the smallest number of sets in $\mathcal{I}$ with union not in $\mathcal{I}$. Formally, $\text{add}(\mathcal{I}) := \min\{|A| : A \subseteq I \land \bigcup A \notin I\}.$

    \item The \emph{covering} number of $\mathcal{I}$, $\mathrm{cov}(\mathcal{I})$, is the smallest number of sets in $\mathcal{I}$ with union $X$. Formally, $\text{cov}(\mathcal{I}) := \min\{|A| : A \subseteq I \land \bigcup A = X\}.$

    \item The \emph{uniformity} of $\mathcal{I}$, $\mathrm{non}(\mathcal{I})$, is the smallest cardinality of any subset of $X$ not in $\mathcal{I}$. Formally, $\text{non}(\mathcal{I}) := \min\{|A| : A \subseteq X \land A \notin I\}.$

    \item The \emph{cofinality} of $\mathcal{I}$, $\mathrm{cof}(\mathcal{I})$, is the smallest cardinality of a $\mathbf{B}\subseteq\mathcal{I}$ such that each element of $\mathcal{I}$ is a subset of an element of $\mathbf{B}$. Such a $\mathbf{B}$ is called a \emph{basis} for $\mathcal{I}$. Formally, $\text{cof}(\mathcal{I}) := \min\{|A| : A \subseteq I \land (\forall B \in I)(\exists A \in \mathcal{A})(B \subseteq A)\}.$
\end{enumerate}
\end{definition}

Define, $\mathrm{Cof}(\mathcal{I}):=(\mathcal{I},\mathcal{I},\subseteq)$ and $\mathrm{Cov}(\mathcal{I}):=(X,\mathcal{I},\in)$. Then, $\norm{\mathrm{Cov}(\mathcal{I})}=\mathrm{cov}(\mathcal{I})$, $\norm{\mathrm{Cov}(\mathcal{I})^{\perp}}=\mathrm{non}(\mathcal{I})$, $\norm{\mathrm{Cof}(\mathcal{I})}=\mathrm{cof}(\mathcal{I})$, and $\norm{\mathrm{Cof}(\mathcal{I})^{\perp}}=\mathrm{add}(\mathcal{I})$.

When the ideal $\mathcal{I}$ is generated by first category (meager) sets, we denote it as $\mathcal{B}$, for "Baire." Likewise, if the ideal is generated by Lebesgue measure zero sets, we denote it as $\mathcal{L}$, for "Lebesgue." Additionally, if we take our underlying set to be different versions of the continuum ($\mathbb{R}$, $\prescript{\omega}{\hphantom{.}}{2}$, $\prescript{\omega}{\hphantom{.}}{\omega}$, etc), we will not notationally distinguish them. Each version of all of the relations admits morphisms in both directions, making them essentially equivalent. With this all in mind, we are now ready to give our first result about cardinal characteristics related to ideals. 
\begin{theorem}\label{ideal cor}
There exist morphisms: \begin{enumerate}
    \item $\varphi_1:\mathrm{cof}(\mathcal{B})\rightarrow \mathrm{non}(\mathcal{B})$.
    \item $\varphi_2:\mathrm{cov}(\mathcal{B})\rightarrow \mathrm{add}(\mathcal{B})$.
    \item $\varphi_3:\mathrm{cof}(\mathcal{L})\rightarrow \mathrm{non}(\mathcal{L})$.
    \item $\varphi_4:\mathrm{cov}(\mathcal{L})\rightarrow \mathrm{add}(\mathcal{L})$.
\end{enumerate}
\end{theorem}
\begin{proof}
We want to construct morphisms: 
$$(\mathcal{I},\mathcal{I},\subseteq)\xrightarrow{\varphi_1}(X,\mathcal{I},\in)\xrightarrow{\varphi_2}(\mathcal{I},\mathcal{I},\not\supset)$$

To construct $\varphi_1$, let $\varphi_{1_{-}}$ map any $ x\in X$ to the set $\{x\}$ and let $\varphi_{1_{+}}$ be the identity. For some $J\in\mathcal{I}$, if $\{x\}\subseteq J$ then $x\in J$. 

To construct $\varphi_2$ first let $A,J\in\mathcal{I}$ be arbitrary. Since $\mathcal{I}$ is proper, we can choose $a\in X\setminus A$ arbitrarily. Define $\varphi_{2_{-}}(A)=a$ and let $\varphi_{2_{+}}$ be the identity. If $a\in J$ then $J\not\subset A$, since $a\not\in A$. 

By Theorem \ref{morphism theorem}, taking the duals of $\varphi_1$ and $\varphi_2$ finishes the proof. 
\end{proof}
\begin{corollary}
\label{ideal relations} $\mathrm{add}(\mathcal{I})\leq \mathrm{cov}(\mathcal{I}),
\mathrm{non}(\mathcal{I})\leq \mathrm{cof}(\mathcal{I}).$
\end{corollary}
Theorem \ref{ideal cor} also implies the existence of several other morphisms but, since our goal is to create Cichoń's diagram, we will take these to be superfluous.
\subsection{Baire Characteristics}\hphantom{.}

Our next objective is to construct morphisms witnessing the inequalities $\mathrm{cof(\mathcal{B}})\geq\mathfrak{d}\geq\mathrm{cov}(\mathcal{B})$ and $\mathrm{non}(\mathcal{B})\geq\mathfrak{b}\geq\mathrm{add}(\mathcal{B}).$ To do so we must first define the concepts of "chopped reals" and "matching." This will allow us to prove two lemmas which, in turn, let us construct alternate relations witnessing $\mathrm{cof}(\mathcal{B})$ and $\mathrm{cov}(\mathcal{B})$. Given these alternate relations, we can construct our desired morphisms.

\begin{definition}\label{chopped real}
A \emph{chopped real} is a pair $(f, \Pi)$, where $f \in 2^\omega$ and $\Pi$ is an interval partition of $\omega$. We write "$\mathrm{CR}$" for the set $\prescript{\omega}{\hphantom{.}}{2} \times \mathrm{IP}$ of chopped reals. We say a real $y \in 2^\omega$ \emph{matches} a chopped real $(x, \Pi)$ if $x \vert I = y \vert I$ for infinitely many $I \in \Pi$. 
\end{definition}

\begin{lemma}\label{Meager lemma}
$M\subseteq\prescript{\omega}{\hphantom{.}}{2}$ is meager if and only if there is a chopped real that no member of $M$ matches.
\end{lemma}
\begin{proof}
For $n,k\in\omega$, the set of all reals matching a chopped real $(f,\Pi)$ is: $$\mathrm{Match}(f,\Pi)=\bigcap_{k=0}^\infty\bigcup_{n\geq k}\{h\in \prescript{\omega}{\hphantom{.}}{2}:f\vert I_n=h\vert I_n\text{ for $I_n\in\Pi$}\}.$$

In the product topology on $\prescript{\omega}{\hphantom{.}}{2}$, basic open sets are families $\mathcal{O}\subseteq \prescript{\omega}{\hphantom{.}}{2}$ which are fixed for finitely many coordinates. Therefore, each set $\{ h \in \prescript{\omega}{\hphantom{.}}{2} : h \vert I_n = f \vert I_n \}$ is open. To show that $D:={\bigcup_{n \geq k} \{ h \in \prescript{\omega}{\hphantom{.}}{2}}: h \vert I_n = f \vert I_n \}$ is dense, consider an open set $\mathcal{O}_X\subset\prescript{\omega}{\hphantom{.}}{2}$ defined by fixing values on some finite set $X\subset\omega$. For some $j_0\in\omega$, let $x\in I_{j_0}$ be the maximum element of $X$. Let $h\in\prescript{\omega}{\hphantom{.}}{2}$ be a function fulfilling the open set conditions on $X$ and require that $h|I_{j_1}=f|I_{j_1}$, for some $j_1>\max\{j_0,k\}$. By construction, $h\in D\cap\mathcal{O}_X$. Since $\mathcal{O}_X$ was arbitrary, this shows $D$ is dense. Thus, $\mathrm{Match}(f,\Pi)$ is a countable intersection of dense open sets. By the Baire Category Theorem, $\prescript{\omega}{\hphantom{.}}{2}\setminus\text{Match}{(f,\Pi)}$ is meager and so $M\subseteq\prescript{\omega}{\hphantom{.}}{2}\setminus\text{Match}{(f,\Pi)}$ is meager. 

A nowhere dense set $Y\subseteq\prescript{\omega}{\hphantom{.}}{2}$ is such that for every finite binary sequence $u\in\prescript{<\omega}{\hphantom{.}}{2}$, there exists a finite extension $u^\frown w$, such that no $y\in Y$ further extends $w$. To prove the other direction, first assume that $M\subseteq\prescript{\omega}{\hphantom{.}}{2}$ is meager. Let $F_n$ be a countable sequence of nowhere-dense sets which cover $M$. For each $i\in\omega$, it can be assumed that $F_i\subseteq F_{i+1}$. Our goal is to construct a chopped real $(f,\Pi:=\{I_n:n\in\omega\})$, such that for each $n\in\omega$ and $t\in F_n$, $f\vert I_n\neq t\vert I_n$. If we can do this, since any $g\in M$ is such that $g\in F_i$ for some $i\in\omega$ and $F_i\subseteq F_{i+1}\subseteq\dots$, we would have that $g\vert I_n \neq f \vert I_n$ for all $n\geq i$, i.e. $g\not\in\mathrm{Match}(f,\Pi)$.

To construct such a chopped real we use recursion. Assume that $I_k$ and $f\vert I_k$ are defined for each $k\leq n-1$. If
$r$ is the right endpoint of $I_{n-1}$, we can let $A:=\{a_i\in\prescript{r}{\hphantom{.}}{2}:i\in 2^r\}$ enumerate all $r+1$ length binary sequences. We want to construct $I_n$ using $2^r$ non-overlapping sub-intervals $J_i$, so that $I_n=\cup_{i\in 2^r}J_i$. Since $F_n$ is nowhere dense, there is some function $h_0\in\prescript{\omega}{\hphantom{.}}{2}$ which extends $a_0$ onto $J_0$, so that none of $F_n$ further extends $a_0^\frown h_0$. Define $f\vert J_0$ as $h_0\vert J_0$. Likewise, there is some function $h_1\in\prescript{\omega}{\hphantom{.}}{2}$ which extends $a_1^\frown h_0$ onto $J_1$, so that none of $F_n$ further extends $a_1^\frown h_0^\frown h_1$. Define $f\vert J_1$ as $h_1\vert J_1$, and so on. By recursively constructing in this way, after $2^r$ steps we will have defined both $I_n$ and $f\vert I_n$. Let $n\in\omega$ and $g\in F_n$ be arbitrary. By construction $f\vert{[0,r)}=a_i$, for some $i\in\omega$ and so $f\vert J_i\neq g\vert J_i$. Since $J_i\subset I_n$, $f\vert I_n\neq g\vert I_n$.
\end{proof}

\begin{lemma}\label{Match lemma}$\mathrm{Match}(f,\Pi)\subseteq\mathrm{Match}(g,\Pi')$ if and only if for all but finitely many intervals $I\in\Pi$ there exists an interval $J\in\Pi'$ such that $J\subseteq I$ and $g\vert J=f\vert J$.
\end{lemma}
\begin{proof}
First assume that $\mathrm{Match}(f,\Pi)\subseteq\mathrm{Match}(g,\Pi')$. For contradiction, suppose that there are infinitely many $I\in\Pi$ such that for each $J\in\Pi'$, $J\subseteq I$ implies $f\vert J\neq g\vert J$. Let $Q$ be an infinite set of such $I\in\Pi$ and let $h\in\prescript{\omega}{\hphantom{.}}{2}$ agree with $f$ on each $I\in Q$. Moreover, for each $i\in\omega\setminus\cup Q$, define $h(i)=1-g(i)$.

Since $h$ agrees with $f$ on all of $\cup Q$, it has to agree with $f$ on infinitely many $I\in\Pi$. Thus,  $h\in\mathrm{Match}(f,\Pi)$ and by assumption $h\in\mathrm{Match}(g,\Pi')$. Define the set $Y:=\{J\in\Pi':h\vert J=g\vert J\}$, which is infinite because $h\in\mathrm{Match}(g,\Pi')$. Notice that each $J\in Y$ is a subset of $\cup Q$. If $J\subseteq\omega\setminus\cup Q$, by construction $h\vert J\neq g\vert J$. But, this contradicts the fact that $J\in Y$. Thus, for each $J\in Y$ there exists an $I\in \Pi$ such that $J\subseteq I$. So given a $J\in Y$, by assumption $f\vert J\neq g\vert J$ but, since $J\subseteq\cup Q$, $f\vert J=h\vert J=g\vert J$.

To prove the reverse direction assume that  $\mathrm{Match}(f,\Pi)\not\subseteq\mathrm{Match}(g,\Pi')$. Let $h\in\mathrm{Match}(f,\Pi)\setminus\mathrm{Match}(g,\Pi')$ and suppose that for all but finitely many $I\in\Pi$, there exists a $J\in\Pi'$ such that $J\subseteq I$ and $g\vert J=f\vert J$. Define $X:=\{I\in\Pi:h\vert I=f\vert I\}$ and observe that $X$ is infinite since $h\in\mathrm{Match}(f,\Pi)$. By assumption, for each of the infinitely many $I\in X$, there are infinitely many $J\in\Pi'$ with $J\subseteq I$,  such that $g\vert J=f\vert J$. Consequently, there are infinitely many $J\in\Pi'$ such that $g\vert J=f\vert J=h\vert J$, which implies that $h\in\mathrm{Match}(g,\Pi')$.\hphantom{\cite{monk}}
\end{proof}
Before we utilize the above two lemmas to construct an alternate witness for  $\mathrm{cof(\mathcal{B})}$, we must first define a binary relation for this
witness. To this end, we give the following definition.
\begin{definition}
 A chopped real $(f,\Pi)$ \emph{engulfs} another chopped real $(x',\Pi')$ if, $\mathrm{Match}(f,\Pi)\subseteq\mathrm{Match}(g,\Pi').$  
\end{definition}
 
 In our construction for an alternate witness to $\mathrm{cof(\mathcal{B})}$ will follow naturally to define an alternate witness for $\mathrm{cov}(\mathcal{B})$. So we will construct both in the following lemma.

\begin{lemma}\label{morphas}
There exist morphisms 
\begin{enumerate}
    \item \(\varphi_1: \mathrm{Cof}(\mathcal{B}) \to (\mathrm{CR}, \mathrm{CR}, \text{is engulfed by})\)
    \item \(\varphi_2: (\mathrm{CR}, \mathrm{CR}, \text{is engulfed by}) \to \mathrm{Cof}(\mathcal{B})\)
    \item \(\varphi_3: \mathrm{Cov}(\mathcal{B}) \to (\prescript{\omega}{\hphantom{.}}{2}, \mathrm{CR}, \text{does not match})\)
    \item \(\varphi_4: (\prescript{\omega}{\hphantom{.}}{2}, \mathrm{CR}, \text{does not match})\to \mathrm{Cov}(\mathcal{B})\)
\end{enumerate}

\begin{proof} 
\hphantom{.}
\begin{enumerate}
\item As in Lemma \ref{Meager lemma}, define $\varphi_{1_{-}}:\mathrm{CR}\rightarrow\mathcal{B}$ as $\varphi_{1_{-}}((f,\Pi)):=\prescript{\omega}{\hphantom{.}}{2}\setminus\mathrm{Match}(f,\Pi)$. Define $\varphi_{1_{+}}:\mathcal{B}\rightarrow\mathrm{CR}$ as $\varphi_{1_{+}}(b)=(g,\Pi)_{b}$, where $(g,\Pi)_{b}$ is a chopped real such that no member of $b$ matches it. Such a chopped real exists because of Lemma \ref{Meager lemma}. Assume that $\prescript{\omega}{\hphantom{.}}{2}\setminus\mathrm{Match}(f,\Pi)\subseteq b$ for an arbitrary chopped real $(f,\Pi)$ and $b\in\mathcal{B}$. We want to show that $(f,\Pi)$ is engulfed by $(g,\Pi)_{b}$. If $h\in\prescript{\omega}{\hphantom{.}}{2}\setminus\mathrm{Match}(f,\Pi)$ then $h\in b$, which implies
$h\not\in\mathrm{Match}(g,\Pi)_{b}$. By contraposition, if $h\in\mathrm{Match}(g,\Pi)_{b}$ then $h\in\mathrm{Match}(f,\Pi)$. By definition, $(f,\Pi)$ is engulfed by $(g,\Pi)_{b}$. 

\item As above, define $\varphi_{2_{-}}:\mathcal{B}\rightarrow\mathrm{CR}$ as $\varphi_{2_{-}}(b)=(g,\Pi)_{b}$ for any $b\in\mathcal{B}.$ Define $\varphi_{2_{+}}:\mathrm{CR}\rightarrow\mathcal{B}$ as $\varphi_{2_{+}}((f,\Pi))=\prescript{\omega}{\hphantom{.}}{2}\setminus\mathrm{Match}(f,\Pi)$ for any $(f,\Pi)\in\mathrm{CR}$. If we assume that $(g,\Pi)_{b}$ is engulfed by  $(f,\Pi)$, Lemma \ref{Match lemma} says $\mathrm{Match}(f,\Pi)\subseteq\mathrm{Match}(g,\Pi)_{b}$. If $x\in b$ then $x\not\in\mathrm{Match}(g,\Pi)_{b}$. This implies that $x\not\in\mathrm{Match}(f,\Pi)$, i.e. $b\in\prescript{\omega}{\hphantom{.}}{2}\setminus\mathrm{Match}(f,\Pi)$.

\item Let $\varphi_{3_{-}}:\prescript{\omega}{\hphantom{.}}{2}\rightarrow\prescript{\omega}{\hphantom{.}}{2}$ be the identity and  for any $b\in\mathcal{B}$ let $\varphi_{3_{+}}:\mathcal{B}\rightarrow\mathrm{CR}$ equal $(g,\Pi)_b$. If $f\in b$, then $f$ does not match $(g,\Pi)_b$.

\item Let $\varphi_{4_{-}}:\prescript{\omega}{\hphantom{.}}{2}\rightarrow\prescript{\omega}{\hphantom{.}}{2}$ be the identity and let  $\varphi_{4_{+}}:\mathrm{CR}\rightarrow\mathcal{B}$ be defined as $\varphi_{4_{+}}((f,\Pi))=\prescript{\omega}{\hphantom{.}}{2}\setminus\text{Match}(f,\Pi)$. If $g\in\prescript{\omega}{\hphantom{.}}{2}$ does not match $(f,\Pi)\in\mathrm{CR}$, then $g\in\prescript{\omega}{\hphantom{.}}{2}\setminus\text{Match}(f,\Pi)$.

\end{enumerate}
\end{proof}
\end{lemma}

With the two alternate relations in mind, we are finally ready to prove the inequalities we initially sought.

\begin{theorem}\label{add<b d<cofB}$\mathrm{add}(\mathcal{B})\leq\mathfrak{b}\leq\mathrm{non}(\mathcal{B})$ and $\mathrm{cof}(\mathcal{B})\geq\mathfrak{d}\geq\mathrm{cov}(\mathcal{B})$.
\begin{proof}
We first produce a morphism from $\mathrm{Cof}(\mathcal{B})$ into $\mathfrak{D'}$. For an arbitrary $f\in\prescript{\omega}{\hphantom{.}}{2}$,  define $\varphi_{-}:\mathrm{IP}\rightarrow\mathrm{CR}$ as $\varphi_{-}(\Pi)=(f,\Pi)$. Define $\varphi_{+}:\mathrm{CR}\rightarrow\mathrm{IP}$ as $\varphi_{+}(g,\Pi')=\Pi'$. Assume that $(f,\Pi)$ is engulfed by some $(g,\Pi')\in\mathrm{CR}$. By Lemma \ref{Match lemma}, this implies that $\Pi'$ must dominate $\Pi$.

We now construct a morphism from $\mathfrak{D}$ into $\mathrm{Cov(\mathcal{B)}}$. Considering the duals of this morphism and the above morphism will complete the proof. By taking both $\phi_+$ and $\phi_-$ to be the identity, we obtain a morphism between $\mathfrak{D}$ and the relation $\mathcal{W}=(\prescript{\omega}{\hphantom{.}}{\omega},\prescript{\omega}{\hphantom{.}}{\omega},\text{is eventually different than}).$ For $f,g\in\prescript{\omega}{\hphantom{.}}{\omega}$, $f$ is eventually different than $g$ means that for all but finitely many $n\in\omega$, $f(n)\neq g(n)$. Then, Lemma \ref{morphas} implies that constructing a morphism $\psi:\mathcal{W}\rightarrow(\prescript{\omega}{\hphantom{.}}{2}, \mathrm{CR}, \text{does not match})$ and composing it with the above morphism completes the proof. To this end, let $\psi_-$ be the identity and assume that $f\in\prescript{\omega}{\hphantom{.}}{\omega}$ is eventually different than $g\in\prescript{\omega}{\hphantom{.}}{\omega}$. Let $\Pi_0:=\{\{n\}:n\in\omega\}$, defining $\psi_+(g)=(f,\Pi_0)$ implies that $f$ cannot match $\psi_+(g)$. If it did, by the construction of $\Pi_0$, for all but finitely many $n\in\omega$, $f(n)=g(n)$. Thus, $\psi$ defines a morphism from $\mathcal{W}$ to $(\prescript{\omega}{\hphantom{.}}{2}, \mathrm{CR}, \text{does not match})$.
\end{proof}
\end{theorem}

It has been proven \cite{bartoszynski} that $\norm{\mathcal{W}}=\mathrm{cov}(\mathcal{B})$ and $\norm{\mathcal{W}^\perp}=\mathrm{non}(\mathcal{B})$ . So technically, the identity morphism from $\mathfrak{D}$ into $\mathcal{W}$ completed the latter half of the proof.
\subsection{Lebesgue Characteristics}\hphantom{.}

Now that we have constructed morphisms between each of the category relations and measure relations we defined, it only remains to bridge the gap between measure and category.
\begin{theorem}\label{measure-cat}
 $\mathrm{cov}(\mathcal{B})\leq\mathrm{non}(\mathcal{L})$ and $\mathrm{cov}(\mathcal{L})\leq\mathrm{non}(\mathcal{B})$. 
 \begin{proof}
  Our goal is to find a suitable relation $K\subseteq\prescript{\omega}{\hphantom{.}}{2}\times\prescript{\omega}{\hphantom{.}}{2}$, such that we can construct morphisms $\varphi$, from the  relation $\mathbf{K}:=(\prescript{\omega}{\hphantom{.}}{2},\prescript{\omega}{\hphantom{.}}{2},K)$ into $\mathrm{Cov}(\mathcal{L})$ and $\psi$, from $\mathbf{K}^\perp$ into $\mathrm{Cov}(\mathcal{B}).$ Then $\varphi^\perp\circ\psi$ and $\psi^\perp\circ\varphi$ witness our desired inequalities. To this end, let $\Pi$ be an interval partition of $\omega$ whose $n^{th}$ interval has $n+1$ elements. Take $fKg$ to mean that there are infinitely many $n\in\omega$ such that $f\vert I_n=g\vert I_n$. Moreover, define $K_f:=\{g\in\prescript{\omega}{\hphantom{.}}{2}:fKg\}$ 
  
  Before we construct the morphisms, we will prove $K_f$ is co-meager and measure zero. The claim that $K_f$ is co-meager follows directly from Lemma \ref{Meager lemma}, since none of $\omega\setminus K_f$ matches $(f,\Pi)$. To prove $K_f$ has measure zero, for every $n\in\omega$, define $A_n:=\{g\in\prescript{\omega}{\hphantom{.}}{2}:\exists k\geq n,\:g\vert I_k=f\vert I_k\}$ and notice $K_f=\bigcap^{\infty}_{n=1}A_n$. Since open sets fix a finite number of coordinates, we can place a probability measure on  $\prescript{\omega}{\hphantom{.}}{2}$ based on a given finite sequence. Concretely, if we specify $j\in\omega$ bits of a sequence, the probability that any sequence agrees with this partial assignment is $2^{-j}$. Extending this to infinitely many coordinates gives the full probability measure on $\prescript{\omega}{\hphantom{.}}{2}$. In our case, a single interval $I_k$ has length $k+1$ and the set of all sequences that match $f$ on $I_k$ has measure $2^{-(k+1)}$. Recalling the definition of $A_n$, this implies $\mu(A_n)\leq \sum^{\infty}_{k=n}2^{-(k+1)}=2^{-n}$. By continuity of the measure, $$\mu(K_f)=\mu(\bigcap^{\infty}_{n=1}A_n)\leq\lim_{n\rightarrow\infty}(A_n)\leq\lim_{n\rightarrow\infty}2^{-n}=0$$.

 Let  both $\varphi_-$ and $\psi_-$ be the identity on $\prescript{\omega}{\hphantom{.}}{2}$. If we let $\varphi_+(f)=K_f$ and $\psi_+(f)=\prescript{\omega}{\hphantom{.}}{2}\setminus K_f$, we get our desired morphisms. 
 \end{proof}
\end{theorem}

It only remains for us to prove there exists a morphism witnessing the inequalities $\mathrm{add}(\mathcal{L})\leq \mathrm{add}(\mathcal{B})$ and $\mathrm{cof}(\mathcal{B})\leq \mathrm{cof}(\mathcal{L})$. We start with the following lemma. 

\begin{lemma}\label{cof add lemma}
$\mathrm{cof}(\mathcal{B})=\text{max}\{\mathrm{non}(\mathcal{B}),\mathfrak{d}\}$ and $\mathrm{add}(\mathcal{B})=\text{min}\{\mathrm{cov}(\mathcal{B}),\mathfrak{b}\}$
\begin{proof}
$\mathrm{cof}(\mathcal{B})\geq\text{max}\{\mathrm{non}(\mathcal{B}),\mathfrak{d}\}$ and $\mathrm{add}(\mathcal{B})\leq\text{min}\{\mathrm{cov}(\mathcal{B}),\mathfrak{b}\}$ is already implied by Theorem \ref{ideal cor} and Theorem \ref{add<b d<cofB}. So if we can exhibit a morphism $$\varphi:(\mathrm{CR},\prescript{\omega}{\hphantom{.}}{2}, \text{ matches)};\mathfrak{D'}\rightarrow(\mathrm{CR},\mathrm{CR},\text{is engulfed by}),$$ by Theorem \ref{extra relation operations} and Lemma \ref{morphas}, we will have proved the lemma. Since $\mathrm{CR}:=\prescript{\omega}{\hphantom{.}}{2}\times\mathrm{IP}$, we can allow $\varphi_+$ to be the identity. To define $\varphi_-$ we require two maps $f:\mathrm{CR}\rightarrow\mathrm{CR}$ and $g:\mathrm{CR} \rightarrow \prescript{\prescript{\omega}{\hphantom{.}}{2}}{\hphantom{.}}{\mathrm{IP}}$. We let $f$ be the identity and define $g$
as follows. For each $(k,\Pi_k)\in\mathrm{CR}$
and each \(h\in\prescript{\omega}{\hphantom{.}}{2}\), if $h$ does not match $(k,\Pi_k)$, let $g((k,\Pi_k))(h)\in\mathrm{IP}$ be arbitrary.  If $h$ does match $(k,\Pi_k)$, let $g((k,\Pi_k))(h)$ be such that each of its blocks contains exactly one of the infinitely many intervals where $k$ and $h$ agree.

Define $\varphi_{-}((x,\Pi_x)):=((x,\Pi_x),g)$ and let $(y,\Pi_y)$ be arbitrary. We require that if $y\text{ matches }(x,\Pi_x)$ and $g((x,\Pi_x))(y)$ is dominated by $\Pi_y$, then $(x,\Pi_x)$ is engulfed by $(y,\Pi_y)$. This implication trivially holds if $y$ does not match $(x,\Pi_x)$, so we assume it does. By construction, within each block $I_g\in g(x,\Pi_x)(y)$ there exists an interval $I_x\subseteq I_g$ for which $x\vert I_x=y\vert I_x$. Since $\Pi_y$ dominates $g(x,\Pi_x)(y)$, for all but finitely many $I_g\in g(x,\Pi_x)(y)$ there exists some $I_y\in\Pi_y$ such that $I_g\subseteq I_y$. Thus, for all of these $I_g$, we have $I_x\subseteq I_g\subseteq I_y$, where $x|I_x=y|I_x$. By Lemma \ref{Match lemma}, $(x,\Pi_x)$ is engulfed by $(y,\Pi_y).$
\end{proof}
\end{lemma}

Before we construct our desired morphism, we need to give an alternate characterization of $\mathrm{add}(\mathcal{L})$, first discovered in \cite{Bartoszynski2}. To this end, we give the following definition.

\begin{definition}\label{slalom}
A \emph{slalom} is a function $S:\omega\rightarrow\mathcal{P}(\omega)$ such that for each $n\in\omega$, $S(n)\subset\omega$ has cardinality $n$. We say that a function $f\in\prescript{\omega}{\hphantom{.}}{\omega}$, \emph{goes through} a slalom $S$ if, for all but finitely many $n\in\omega$, $f(n)\in S(n)$. We denote $\mathbb{S}$ as the set of all slaloms. 
\end{definition}

\begin{proposition}\label{alt addL}
$\norm{\mathrm{add}(\mathcal{L})}=(\mathbb{S},\prescript{\omega}{\hphantom{.}}{\omega}, \text{does not go through})$. 
\end{proposition}

The proof of this theorem can be found in \cite{Bartoszynski2}. There it is proven that $\mathrm{add}(\mathcal{L})$ is the least cardinality of any family of functions such that there is no single slalom through which all of the members of it go. This is equivalent to the above norm condition.

\begin{theorem}
 $\mathrm{add}(\mathcal{L})\leq \mathrm{add}(\mathcal{B})$ and $\mathrm{cof}(\mathcal{B})\leq \mathrm{cof}(\mathcal{L})$.   
 \begin{proof}
Given Lemma \ref{cof add lemma}, it suffices to exhibit morphisms witnessing $\mathfrak{b}\geq \mathrm{add}(\mathcal{L})$ and $\mathrm{cov}(\mathcal{B})\geq \mathrm{add}(\mathcal{L})$. Since $\mathrm{add}(\mathcal{B})=\text{min}\{\mathrm{cov}(\mathcal{B}),\mathfrak{b}\}$, by the universal property of co-products, this would show that there exists a unique morphism witnessing  $\mathrm{add}(\mathcal{L})\leq \mathrm{add}(\mathcal{B})$. The dual of this morphism will witness $\mathrm{cof}(\mathcal{B})\leq \mathrm{cof}(\mathcal{L})$. 

We first construct a morphism $\varphi_1:(\prescript{\omega}{\hphantom{.}}{\omega},\prescript{\omega}{\hphantom{.}}{\omega},\not>^*)\rightarrow(\mathbb{S},\prescript{\omega}{\hphantom{.}}{\omega}, \text{does not go through})$. Let $\varphi_{1_-}:\mathbb{S}\rightarrow\prescript{\omega}{\hphantom{.}}{\omega}$ map a slalom $S$ to the function $f_S$, where $f_S(n)=\text{max}\{S(n)\}+1$, for every $n\in\omega$. For some $g\in\prescript{\omega}{\hphantom{.}}{\omega}$, if we suppose that $g\not<^*f_s$, then there are infinitely many $n\in\omega$ such that $g(n)\geq f_s(n)$. By construction of $f_s$, there are infinitely many $n\in\omega$ such that,  $g(n)\not\in S(n)$. Thus, defining $\varphi_+$ as the identity suffices. 

The remark at the end of Theorem \ref{add<b d<cofB}, states that $\mathrm{Cov}(\mathcal{B})$ is equivalent to the relation $\mathcal{W}=(\prescript{\omega}{\hphantom{.}}{\omega},\prescript{\omega}{\hphantom{.}}{\omega},\text{is eventually different than}).$ The most straightforward way to produce a morphism $\varphi_2:\mathrm{Cov}(\mathcal{B})\rightarrow\mathrm{Add}(\mathcal{L})$ is through $\mathcal{W}$. But, because we have not justified the equality of $\mathrm{Cov}(\mathcal{B})$ with $\mathcal{W}$, we refer to the proof given in chapter 5 of \cite{blass} for this morphism. Although the morphism is not given explicitly in the proof, it is not difficult to implicitly draw it out.
 \end{proof}
\end{theorem}

\section{More Cardinal Characteristics}

Although Cichoń's diagram is a useful classification tool for cardinal characteristics related to category and measure, many more cardinal characteristics exist beyond the ones mentioned in the previous section. For many of these cardinal characteristics $\kappa$, apart from the trivial relation \((\kappa,\kappa,=)\), it is often challenging to define a relation that is easy to work with. The purpose of the \(\mathbb{GT}\) category is to simplify the study of cardinal characteristics of the continuum. Therefore, if incorporating a certain cardinal characteristic or inequality into this framework complicates the proofs without yielding any new insights, then its inclusion is unnecessary. In this section, we shift our attention towards a subset of cardinal characteristics, beyond the ones mentioned in Cichoń's diagram, which are easily amenable to the $\mathbb{GT}$ framework.

\subsection{Splitting and Reaping}\hphantom{.}

As defined in Example \ref{splitting example}, the first two cardinal characteristics we will discuss are $\mathfrak{s}$ and $\mathfrak{r}$.

When we first introduced $\mathfrak{b}$ and $\mathfrak{d}$, we sought to construct a morphism from a relation whose norm is $\omega_1$ into a relation whose norm is $\mathfrak{b}$. However, as indicated by Proposition \ref{noMorph}, such a morphism would necessitate a very non-standard witness for $\omega_1$, so we instead opted for a direct proof. In a nearly identical way, we now show that attempting to build a morphism from a relation with norm $\omega_1$ into a relation with norm $\mathfrak{s}$ necessitates a witness for $\omega_1$ which is just as non-standard.

\begin{proposition}\label{noMorph2}
Let $\mathbf{A}=(A_-,A_+,A)$ be a relation such that $\norm{\mathbf{A}}>1$. If $\varphi:([\omega]^\omega,\mathcal{P}(\omega),\text{split by})\rightarrow\mathbf{A}$ is a morphism, then $|A_-|\geq\mathfrak{r}$.
\end{proposition}
\begin{proof}
For contradiction, suppose $|A_-|<\mathfrak{r}.$ Define $\mathcal{R}:=\{\varphi_-(a):a\in A_-\}\subseteq[\omega]^\omega$. Since $\mathcal{R}$ is not a reaping family, let $x\in[\omega]^\omega$ split every element of $\mathcal{R}$. By the definition of a morphism, for every $a\in A_-$, $aA\varphi_+(x)$. This implies $\norm{\mathbf{A}}=1$. By contradiction, $\mathcal{R}$ is a reaping family. But this is impossible since $|A_-|<\mathfrak{r}$. 
\end{proof}

This proof is nearly identical to Proposition \ref{noMorph} and for similar reasons as with $\mathfrak{b}$, we instead opt for a direct proof to show that $\mathfrak{\omega_1}\leq\mathfrak{s}$. Also, as in Proposition \ref{noMorph}, we can just as easily prove Proposition \ref{noMorph2} by considering the dual morphism $\varphi^\perp$. 

\begin{theorem}
$\mathfrak{\omega_1}\leq\mathfrak{s}$.
\end{theorem}
\begin{proof}
Let $\mathcal{S}:=\{Y_n:n\in\omega\}\subseteq[\omega]^\omega$, be a countable family of infinite sets. We will show that $\mathcal{S}$ is not a splitting family. Define $X_0:=Y_0$ and let $x_0\in X_0$. For each $i\in\{0,\dots,n\}$, assume that $X_0\supseteq X_1\supseteq\dots\supseteq X_n$ are infinite sets such that $X_i\subseteq Y_i$ and $x_0,\dots,x_n$ are distinct elements with each $x_i\in X_i$. If $X_n\cap Y_{n+1}$ is infinite let $X_{n+1}:=X_n\cap Y_{n+1}$, otherwise let $X_{n+1}:=X_n\setminus Y_{n+1}$. In either case choose $x_{n+1}\in X_{n+1}$ distinct from $x_0,\dots,x_n$. By induction, the set $Z:=\{x_n:n\in\omega\}$ is infinite. Moreover, for every $k\in\omega$ the set $Z\setminus\{x_0,\dots,x_{k-1}\}$ is a subset of $X_k$ and either  $X_k\subseteq Y_k$ or $X_k\cap Y_k=\emptyset$. Thus for each $k\in\omega$, either $Z\setminus Y_k$ or $Z\cap Y_k$ is finite.
\end{proof}

The cardinal characteristics $\mathfrak{s}$ and $\mathfrak{r}$ share many similarities with $\mathfrak{b}$ and $\mathfrak{d}$, respectively. In fact, it is consistent with ZFC that $\mathfrak{d} = \mathfrak{r}$ and $\mathfrak{b} = \mathfrak{s}$ \cite{halbeisen}. The following morphism succinctly produces the known inequalities between $\mathfrak{s},\mathfrak{r}$ and $\mathfrak{b},\mathfrak{d}$.

\begin{theorem} $\mathfrak{s}\leq\mathfrak{d}$\label{s<d, b<r} and $\mathfrak{b}\leq\mathfrak{r}$
\begin{proof}
We seek a morphism $\varphi:(\prescript{\omega}{\hphantom{.}}{\omega},\prescript{\omega}{\hphantom{.}}{\omega},<^*)\rightarrow([\omega]^\omega,\mathcal{P}(\omega),\text{split by})$. 

Define $\varphi_-:[\omega]^\omega\rightarrow\prescript{\omega}{\hphantom{.}}{\omega}$ as $\varphi_-(x)=f_x$, where $f_x$ is the unique strictly increasing bijection from $\omega$ to $x$. 

Let $\varphi_+:\prescript{\omega}{\hphantom{.}}{\omega}\rightarrow\mathcal{P}(\omega)$ be defined as $\varphi_+(f)=\sigma_f$, where: 
$$\sigma_f:=\bigcup\{[f^{2n}(0),f^{2n+1}(0)):n\in\omega\}.$$
We assume, without loss of generality, that  $f$ is increasing. If $f_x<^*f$, let $n_0$ be the point past which $f$ dominates $f_x$. For every $k>n_0$ we get that :
$$f^k(0)\leq f_x(f^k(0))<f(f^k(0))=f^{k+1}(0).$$
This implies that $f_x(f^k(0))\in\sigma_f$ if and only if $k$ is even. So when $k$ is even $x\cap\sigma_f$ is infinite. Otherwise, if $k$ is odd, $f_x(f^k(0))\not\in\sigma_f$ and so $\sigma_f\setminus x$ is infinite. By definition, $x$ splits $\sigma_f$. 
\end{proof}
\end{theorem}

We can also relate $\mathfrak{s}$ and $\mathfrak{r}$ to category and measure.

\begin{theorem}\label{split/cichon}
$\mathfrak{s}\leq\mathrm{non}(\mathcal{B}),\,\mathrm{non}(\mathcal{L})
\;\:\text{and}\;\:
\mathfrak{r}\geq \mathrm{cov}(\mathcal{B}),\,\mathrm{cov}(\mathcal{L}).$
\end{theorem}
\begin{proof}
For any $A\in[\omega]^\omega$ let
$S_A :=\{X\in[\omega]^\omega : X\text{ does \emph{not} split }A\}$. Before moving forward, notice that the standard product topology on $[\omega]^\omega$ can be aligned with the standard product topology on $\prescript{\omega}{\hphantom{.}}{2}$. This is done by associating each $X\in[\omega]^\omega$ with the set $\{x(n):n\in\omega\}$, where $x:\omega\rightarrow\{0,1\}$ is the characteristic function for $X$. With this in mind, we claim that $S_A$ is meager and of measure zero. To see this notice that $$S_A=\bigcup_{n=0}^\infty\{X\in[\omega]^\omega:|X\cap A|\leq n\}\cup\bigcup_{n=0}^\infty\{X\in[\omega]^\omega:|A\setminus X|\leq n\}.$$ For every $n\in\omega$, both $\{X\in[\omega]^\omega:|X\cap A|\leq n\}$ and $\{X\in[\omega]^\omega:|A\setminus X|\leq n\}$ are meager and measure zero, which implies $S_A$ is meager and measure zero.

To see why $\{X\in[\omega]^\omega:|X\cap A|\leq n\}$ is measure zero, consider the event that any $a\in A$ is also in $X$. Since each such event occurs with probability $\frac{1}{2}$ independently, and $A$ is infinite, the sum of these probabilities diverges. By the second Borel–Cantelli lemma \cite{Borel}, it follows that almost every $X$ will have infinitely many elements from $A$. Thus, the probability that $X$ contains only finitely many elements from $A$ is zero. We can make a similar argument in the case of $\{X\in[\omega]^\omega:|A\setminus X|\leq n\}$.

To show that $\{\,X : |X \cap A|\le n\}$ is meager, observe that since $A$ is infinite, we can choose additional coordinates in $A$ (beyond the finitely many fixed ones) and force those coordinates to be $1$ in $X$. This ensures $X$ meets $A$ in more than $n$ points. Hence, we can find a smaller open set disjoint from $\{\,X : |X \cap A|\le n\}$, showing that this set is nowhere dense. A countable union of nowhere dense sets is meager, so $\{\,X : |X \cap A|\le n\}$ is meager. We can again make a similar argument in the case of $\{X\in[\omega]^\omega:|A\setminus X|\leq n\}$.

Let $\varphi_-$ be the identity and let $\varphi_+(A)=S_A$. Then $\varphi$ can be used as a morphism from $\mathfrak{R}$ into $\mathrm{Cov}(\mathcal{I})$ and from $\mathfrak{R}$ into $\mathrm{Cov}(\mathcal{B})$.
\end{proof}

\subsection{Ramsey-like Characteristics}\hphantom{.}

Ramsey's theorem states that for any $n,r\in\omega$, for any $X\in[\omega]^\omega$, and any coloring $\pi:[X]^n\rightarrow r$, there exists an infinite subset $H\in[X]^\omega$ such that $H$ is homogeneous for $\pi$. As defined in Example \ref{coloring example}, the homogeneity number, denoted $\mathfrak{hom}_n$, is the Ramsey-theoretic analog of $\mathfrak{s}$. The partition number, denoted $\mathfrak{par}_n$, is the Ramsey-theoretic analog of $\mathfrak{r}$.

Our first task will be to prove that $\omega_1\leq\mathfrak{par}_2$. We opt for a direct proof of this fact for similar reasons as with $\mathfrak{s}$ and $\mathfrak{b}$.

\begin{theorem}\label{w<par}
$\omega_1\leq\mathfrak{par}_2.$
\end{theorem}
\begin{proof}
We need to prove that for any countable family of $2$-colorings, there exists an infinite set $H\in[\omega]^\omega$ which is almost homogeneous for all of it. Instead, we prove a stronger analog of this statement. In particular, we show that for any family of colorings $\pi_k:[\omega]^{n_k}\rightarrow m_k$, where $\{n_k:k\in\omega\}$ and $\{m_k:k\in\omega\}$ are countable sets of integers, there is an $H\in[\omega]^\omega$ which is almost homogeneous for the entire family.

We go by induction. Let $A_0\in[\omega]^\omega$ be an infinite set homogeneous for $\pi_0$; $A_0$ exists by Ramsey's theorem. Assume that $A_0,...,A_k$ are defined up to some $k\in\omega$ and that $\pi_0,...,\pi_k$ are homogeneous for each of them respectively. Let $m_k:=\text{min}({A_k})$ and $B_k:=A_k\setminus m_k$. Another application of Ramsey's theorem ensures there is an infinite subset of $B_k$ which is homogeneous for $\pi_{k+1}$, call it $A_{k+1}$. Then, the infinite set $A:=\{m_i:i\in\omega\}$ is almost homogeneous for each coloring in the family.
\end{proof}

The next two propositions seek to solidify the comparisons among $\mathfrak{d,r},\mathfrak{hom}_n$ and $\mathfrak{b,s},\mathfrak{par}_n$. After this, we prove that there exist morphisms witnessing the inequalities $\mathfrak{hom}_2\geq\text{max}\{\mathfrak{r},\mathfrak{d}\}$ and  $\mathfrak{par}_2\leq\text{min}\{\mathfrak{b},\mathfrak{s}\}$. The result actually holds for every $n\geq 2$, but we will only consider when $n=2$. This is because (1) the proofs for arbitrary $n$ are practically identical and (2) the following proposition establishes that the case of $n=2$ suffices to prove all relevant inequalities for $n>2$. 

\begin{proposition}\label{homma}
 For integers $n\geq m$, $\mathfrak{hom_n}\geq\mathfrak{hom_m}$ and $\mathfrak{par_n}\leq\mathfrak{par_m}$.   
 \begin{proof}
It suffices to produce a morphism $\varphi:(P_n,[\omega]^\omega,H)\rightarrow(P_m,[\omega]^\omega,H)$.

Let $\varphi_-$ send any coloring $\pi_m:[\omega]^m\rightarrow 2$ to the coloring $\pi_n:[\omega]^n\rightarrow 2$, defined as $\pi_n\{x_1,...,x_n\}:=\pi_m\{x_1,...,x_m\}$. Let $\varphi_+$ be the identity. If $\pi_n$ is almost homogeneous for some $A\in[\omega]^\omega$, restricting to $m\leq n$ coordinates remains constant after removing a finite set. 
 \end{proof}
\end{proposition}

\begin{proposition}\label{r=hom par=s}
$\mathfrak{hom}_1=\mathfrak{r}$ and $\mathfrak{par}_1=\mathfrak{s}$. 
\begin{proof}
To prove this theorem we exhibit morphisms to and from $\mathfrak{R}$ and $\mathfrak{Hom}_1$.

First we produce a morphism $\varphi_{1}:(\mathcal{P}(\omega),[\omega]^\omega,\text{does not split})\rightarrow(P_1,[\omega]^\omega,H)$. Let $\varphi_{1_-}(\pi)=Z$, the set of all $n\in\omega$ for which $\pi(n)=0$. Let $\varphi_{1_+}$ be the identity. If $Z$ does not split some $A\in[\omega]^\omega$, by definition $Z\cap A$ or $A\setminus Z$ is finite.  If $Z\cap A$ is finite then $\pi\vert A$ is almost always $1$. If $A\setminus Z$ is finite then $\pi\vert A$ is almost always $0$. In any instance, $\pi$ is almost homogeneous for $A$.

For the reverse morphism let $\varphi_{2_-}$ send any $A\in[\omega]^\omega$ to its characteristic coloring. Concretely, let $\varphi_{2_-}(A)=\text{char}_A$ where \[
\text{char}_A(x) := 
\begin{cases}
0 & \text{if }  x\in A\\
1 & \text{if } x\not\in A
\end{cases}
\]
Let $\varphi_{2_+}$ be the identity. If $\text{char}_A$ is almost homogeneous for some $B\in[\omega]^\omega$, then the restriction of $\text{char}_A$ to $B$ is either almost always $0$ or $1$, but not both. For contradiction, assume that $A$ splits $B$. Since $A\cap B$ is infinite, $\text{char}_A\vert B$ is almost  homogeneous for $0$. But since $B\setminus A$ is infinite, $\text{char}_A\vert B$ is almost homogeneous for $1$. By contradiction,  $A$ does not split $B$. 
\end{proof}
\end{proposition}

\begin{theorem}\label{max(r,d)=hom}
    $\text{max}\{\mathfrak{r},\mathfrak{d}\}\leq\mathfrak{hom}_2$ and $\mathfrak{par}_2\leq\text{min}\{\mathfrak{b},\mathfrak{s}\}$
\begin{proof}

It suffices to construct a morphism from $\varphi:(P_2,[\omega]^\omega,H)\rightarrow(\prescript{\omega}{\hphantom{.}}{\omega},\prescript{\omega}{\hphantom{.}}{\omega},<^*)$. This is because composing the morphisms given in Proposition \ref{homma} and Proposition \ref{r=hom par=s} gives a morphism witnessing $\mathfrak{r}\leq\mathfrak{hom}_2$. It is also
possible to do this directly \cite{halbeisen}.

For a strictly increasing function $f\in\prescript{\omega}{\hphantom{.}}{\omega}$, let $\varphi_{-}:\prescript{\omega}{\hphantom{.}}{\omega}\rightarrow P_2$ be defined as $\varphi_{-}(f)=\pi_f$ where, $$\pi_f(\{n,m\})=0\iff\exists k\in\omega:(f(2k)\leq n, m<f(2k+2)).$$

For any $x\in[\omega]^\omega$, define $\varphi_{+}:[\omega]^\omega\rightarrow\prescript{\omega}{\hphantom{.}}{\omega}$ as $\varphi_{+}(x)=f_x$, where $f_x$ is the unique strictly increasing bijection from $\omega$ into $x$. Suppose $x$ is almost homogeneous for $\pi_f$. If $\pi_f(y)=0$ for almost every $y\in[x]^2$, then every pair must lie within the same interval. This is impossible as $x$ is infinite and such an interval would be finite. So it must be that $\pi_f(y)=1$ for almost every $y\in[x]^2$. In the words, for almost every pair $\{n,m\}\in[x]^2$, there is no interval that both $n$ and $m$ are elements of. Notice that the intervals $[f(2k),f(2k+2))$ form a partition of $\omega$ and the least interval that $f_x(0)$ could lie in is $[f(0),f(2))$. So even in the worst case, $f_x(0)\geq f(0)$. Likewise, the least interval that $f_x(1)$ could lie in is $[f(2),f(4))$ and this still implies that $f_x(1)> f(1)$. Continuing in this way, we see that past $n=0$, $f(n)<f_x(n)$. Thus, $f$ is dominated by $f_x$.
 \end{proof}
 \end{theorem}

The above inequality can be strengthened into an equality if we consider $\mathfrak{r}_\sigma$ instead of $\mathfrak{r}$. Particularly, for $n\geq2$, it can be shown that there is a morphism witnessing $\mathfrak{hom}_n\leq\text{max}\{\mathfrak{r_\sigma},\mathfrak{d}\}$ and $\mathfrak{par}_n\geq\text{min}\{\mathfrak{b},\mathfrak{s_\sigma}\}$. Since $\mathfrak{s_\sigma}\geq\mathfrak{s}$, the above theorem implies that $\text{max}\{\mathfrak{r},\mathfrak{d}\}\leq\mathfrak{hom}_n\leq\text{max}\{\mathfrak{r_\sigma},\mathfrak{d}\}$ and $\mathfrak{par}_n=\text{min}\{\mathfrak{b},\mathfrak{s}\}$. Notice that in requiring $\mathfrak{r}_\sigma$ instead of $\mathfrak{r}$ we fail to get a perfectly dual inequality for $\mathfrak{hom}_n$. Fortunately, it still is the case that $\mathfrak{hom}_n=\text{max}\{\mathfrak{r_\sigma},\mathfrak{d}\}$. We just have to construct an extra morphism to prove it. If we attempted to construct the morphism witnessing $\mathfrak{hom}_n\leq\text{max}\{\mathfrak{r_\sigma},\mathfrak{d}\}$, Definiton \ref{alternate definitions} would tell us it must go from $\mathfrak{R}_\sigma;(\mathfrak{R}\land\mathfrak{D} )$ into $\mathfrak{hom}_2$. This would amount to producing a morphism $$\varphi:(\prescript{\omega}{\hphantom{.}}{\mathcal{P}(\omega)\times\prescript{[\omega]^\omega}{\hphantom{.}}{[\omega]^\omega}}\times\prescript{[\omega]^\omega}{\hphantom{.}}{(\prescript{\omega}{\hphantom{.}}{\omega}}),[\omega]^\omega\times[\omega]^\omega\times\prescript{\omega}{\hphantom{.}}{\omega},K)\rightarrow(P_n,[\omega]^\omega,H),$$ where $(a,b,c)K(d,e,f)$ requires three separate conditions to be met. This would be inelegant, to say the least. Alternatively, we provide a direct proof of this result and note that the above morphism can be found implicitly within the argument.

\begin{theorem}\label{parra}
 For $n\geq2$, $\text{max}\{\mathfrak{r_\sigma},\mathfrak{d}\}=\mathfrak{hom}_n$ and $\mathfrak{par}_n=\text{min}\{\mathfrak{b},\mathfrak{s}\}$. 
\begin{proof}\hphantom{.}
We first prove $\text{max}\{\mathfrak{r_\sigma},\mathfrak{d}\}\leq\mathfrak{hom}_2$. To do this, by Theorem \ref{max(r,d)=hom}, it suffices to show $\mathfrak{r}_\sigma \le \mathfrak{hom}_2.$ To this end, we employ the following argument due to Brendle \cite{Brendle1995}.

Let $\mathcal{H}\subseteq[\omega]^\omega$ have $\mathfrak{hom}_2$ property. This means that for any $2$-coloring $\pi:[\omega]^2\to2,$ there exists $H\in\mathcal{H}$ which is almost homogeneous for $\pi.$ 
Given a countable family $(f_n)_{n\in\omega}$ of maps $f_n:\omega\to\{0,1\},$ define for each $x\in\omega$ an infinite binary sequence
$b_x:=(f_0(x),\,f_1(x),\,f_2(x),\,\dots).$
For $x<y$, define a $2$-coloring $L:[\omega]^2\to2$ by
\[
L(\{x,y\}) =
\begin{cases}
0 & \text{if } b_x <_{\mathrm{lex}} b_y\\
1 & \text{else}
\end{cases}
\]
Let $H\in\mathcal{H}$ be almost homogeneous for $L.$ By removing finitely many elements, we obtain an infinite $H'\subseteq H$ on which $L$ is homogeneous in $0$ (a similar argument works for $1$). As $x\in H'$ increases, $f_0(x)$ can only increase. Increasing here means that $f_0(x)$ goes from $0$ to $1$. But once $f_0$ changes from $0$ to $1$, it must remain as $1$. Otherwise, the lexicographical ordering fails. Once $f_0$ stabilizes, $f_1$ increases and then stabilizes, after this $f_2$ increases and stabilizes, and so on. This shows that the family $(f_n)_{n\in\omega}$ is almost homogeneous on $H'$ and thus on $H$.

We now have to prove that $\text{max}\{\mathfrak{r}_\sigma,\mathfrak{d}\}\geq\mathfrak{hom}_2$. To do this we will construct a set $\mathcal{H}$ with cardinality $\text{max}\{\mathfrak{r}_\sigma,\mathfrak{d}\}$ and show that it has $\mathfrak{hom}_2$ property.

Let $\mathcal{D}\subseteq\omega^\omega$ be a dominating family with a cardinality of $\mathfrak{d}$. Let $\mathcal{R}\subseteq[\omega]^\omega$ be a $\sigma$-reaping family with a cardinality of $\mathfrak{r}_\sigma$. Then for each $A\in\mathcal{R}$, we can construct a reaping family $\mathcal{R}_A\subseteq\mathcal{P}(A)$.

To see why this is possible, notice that we can take $\mathcal{R}_A$ to be the set $\{A\cap R:R\in\mathcal{R}\}$. Under this definition, if $\mathcal{R}_A$ is not reaping then we can assume that there is some $X\in[\omega]^\omega$ which splits each of its elements. This would imply that $X$ splits $A\cap R$, for every $R\in\mathcal{R}$. In other words, $X$ would have to split all of $\mathcal{R}$, which is a contradiction.

For every $h\in\mathcal{D}$, $A\in\mathcal{R}$, and $B\in\mathcal{R}_A$, we let $H(h,A,B)\in[B]^\omega$ be such that $h(x)<y$, for any $x<y$ within $H(h,A,B)$. The family $\mathcal{H}$, of each $H(h,A,B)$, has cardinality at most $|\mathfrak{d}|\times|\mathfrak{r}_\sigma|=\text{max}\{\mathfrak{d,r_\sigma}\}$.

Let $\pi$ be an arbitrary $2$-coloring. For each $n\in\omega$ we can define the function $\gamma_n:\omega\rightarrow 2$ as $\gamma_n(x):=\pi\{n,x\}$. Since $\mathcal{R}$ is $\sigma$-reaping, we can find some set $A$ for which each $\gamma_n$ is almost homogeneous on. On $A$, we can say $\gamma_n(x)=j(n)$ for every $x\geq g(n)$. Here, $g(n)$ indicates the point past which $f_n$ is homogeneous and $j(n)$ indicates the color in which $f_n$ is almost homogeneous. Since $\mathcal{R}_A$ is reaping, there is some $B\in\mathcal{R}_A$ such that $j$ is almost constant on it, say $j(n)=i$ for every $n\geq x_0$. Moreover, let $h\in\mathcal{D}$ dominate $g$ past the point $x_1$.

By the construction of $H(h,A,B)$, for $n>m>\text{max}\{x_0,x_1\}$, $g(n)<h(n)<m$. Thus, $\pi\{n,m\}=\gamma_n(m)=j(n)=i$, i.e. $\pi$ is almost constant on $H(h,A,B)$.
\end{proof}
\end{theorem}
\subsection{The Lower Topology}\hphantom{.}

We now introduce the final two cardinal characteristics included in this thesis. To do this, we first define the "lower topology" on $[\omega]^\omega$. In the lower topology, open sets are defined as families $\mathcal{O}\subseteq[\omega]^\omega$ which are closed under almost subsets. Being closed under almost subsets means that if $X\in\mathcal{O}$ and there is a $Y\in[\omega]^\omega$, such that $Y\setminus X$ is finite, then $Y\in\mathcal{O}$. It can be shown that for a family $\mathcal{Y}\subseteq[\omega]^\omega$, being dense is equivalent to every infinite set $A\in[\omega]^\omega$ having a subset $B\subseteq A$ with $B\in\mathcal{Y}$. We let "$\mathrm{DO}$", denote the set of all dense-open sets within the lower topology. 

The first cardinal characteristic we define related to the lower topology is "$\mathfrak{h}$", the "shattering number", also referred to as the "distributivity" number. $\mathfrak{h}$ is the least number of dense-open sets with an empty intersection. Equivalently, we define $\mathfrak{h}:=\norm{\mathfrak{H}}=\norm{([\omega]^\omega,\mathrm{DO},\not\in)}$. 

The second cardinal characteristic related to the lower topology is $\mathfrak{g}$, the "groupwise-dense" number. We say a family $\mathcal{G}\subseteq[\omega]^\omega$ is groupwise dense if it is open within the lower topology and for every interval partition $\Pi\in\mathrm{IP}$, there exists an infinite set of intervals $I\in\Pi$ such that $\cup I\in G$. $\mathfrak{g}$ is the least number of groupwise dense families having empty intersection. Equivalently, we define $\mathfrak{g}:=\norm{\mathfrak{G}}=\norm{([\omega]^\omega,G,\notin)}$, where $G\subseteq[\omega]^\omega$ denotes the set of all groupwise dense subsets. 

Notice that $\mathfrak{g}$ and $\mathfrak{h}$ are not duals of one and other. We will prove that both the duals for $\mathfrak{h}$ and $\mathfrak{g}$ have cardinality $\mathfrak{c}$. To this end, we first introduce the idea of "almost disjoint-ness" and state a proposition relating to it. We say two sets $x,y\in[\omega]^\omega$ are almost disjoint if $x\cap y$ is finite. Moreover, a family $\mathcal{F}\subseteq[\omega]^\omega$ of pairwise almost disjoint infinite sets is called an "almost disjoint family." It should be noted that the least cardinality of a maximal almost disjoint (MAD) family, defines a cardinal characteristic of its own, $\mathfrak{a}$---the "almost disjoint number." It can be shown that $\mathfrak{b}\leq\mathfrak{a}$ but, since defining a relation witnessing this is not obvious to me, we will forego further discussion about $\mathfrak{a}$.   Instead, we continue with the following proposition and lemma.

\begin{proposition}\label{mad family}
There exists a maximal almost disjoint family of cardinality $\mathfrak{c}$.  
\begin{proof}
Let $\prescript{<\omega}{\hphantom{.}}{2}$ denote the set of all finite binary strings. A branch of $\prescript{<\omega}{\hphantom{.}}{2}$ is a function $b:\omega\rightarrow\{0,1\}$. Since $\prescript{\omega}{\hphantom{.}}{\{0,1\}}$ is in bijection with $[0,1]\subset\mathbb{R}$, the set $\mathcal{B},$ of all branches of $\prescript{<\omega}{\hphantom{.}}{2}$ has cardinality $\mathfrak{c}$. For each $b\in\mathcal{B}$, let $D_b:=\{b|n:n\in\omega\}$ be the set of all initial segments of $b$. If $b_0,b_1\in\mathcal{B}$ are distinct then there must be some $k\in\omega$ such that $b_0\vert k\neq b_1\vert k$. This implies $|A_{b_0}\cap A_{b_1}|$ is finite, thus the set $\mathcal{A}:=\{A_b:b\in\mathcal{B}\}$ is almost disjoint with cardinality $\mathfrak{c}$. By Teichmüller's principle, we can extend $\mathcal{A}$ to a maximal almost disjoint family $\mathcal{A'}.$
\end{proof}
\end{proposition}

In addition to being useful for the theorem we want to prove, Proposition \ref{mad family} also shows that $\mathfrak{a}$ is well defined. 

\begin{lemma}\label{MAD lemma}
Let $\mathcal{A}\subseteq[\omega]^\omega$ be a MAD family and let $x\in[\omega]^\omega$. For any distinct elements $A_0, A_1\in\mathcal{A}$ either $x\cap A_0$ or $x\cap A_1$ is finite if and only if there exists an element $A\in\mathcal{A}$ such that $x\setminus A$ is finite.
\begin{proof}
Suppose first that for any two distinct sets $A_0,A_1\in\mathcal{A}$ either $x\cap A_0$ or $x\cap A_1$ is finite. It cannot be the case that every $A\in\mathcal{A}$ has a finite intersection with $x$. If this were the case, then $\mathcal{A}\cup{\{x\}}$ would be a MAD family containing $\mathcal{A}$,  which contradicts maximality. Thus let, $A_i\in\mathcal{A}$ have infinite intersection with $x$. By assumption for every $A\in\mathcal{A}\setminus \{A_i\}$, $A\cap x$ is finite and for each of these $A$, $(x\setminus \{A_i\})\cap A$ is finite. If $x\setminus \{A_i\}$ is infinite then $\mathcal{A}$ would not be MAD, which proves $x\setminus \{A_i\}$ is finite.

For the other direction suppose there is an $A_0\in\mathcal{A}$ such that $x\setminus A_0$ is finite. Since $\mathcal{A}$ is almost disjoint, every $A\in\mathcal{A}\setminus \{A_0\}$ must have a finite intersection with $x$.  
\end{proof}
\end{lemma}

\begin{theorem}\label{contiuumn theorem}
$\norm{\mathfrak{G^\perp}}=\norm{\mathfrak{H^\perp}}=\mathfrak{c}$.
\begin{proof}
We will only prove $\norm{\mathfrak{H}^\perp}=\mathfrak{c}$, as the proof for $\norm{\mathfrak{G}^\perp}$ is almost exactly the same. Since $|[\omega]^\omega|=\mathfrak{c}$, it suffices to prove $\norm{\mathfrak{H}^\perp}\geq\mathfrak{c}$.

 Suppose for contradiction that $\mathcal{Y}:=\{{y_\kappa\in{[\omega]^\omega}:\kappa<\mathfrak{c}}\}$ is the  subset of $[\omega]^\omega$ witnessing $\norm{\mathfrak{H}^\perp}$. Our goal is to construct a dense-open set such that no element of $y$ is in the set. To this end, let $\mathcal{A}$ be a MAD family and define $$\mathcal{A}_D:=\{X\in[\omega]^\omega:\exists A\in\mathcal{A}(X\setminus A\text{ is finite})\}.$$
If $X\in\mathcal{A}_D$ then there is some $A\in\mathcal{A}$ such that $X\setminus A$ is finite. If for some $Y\in[\omega]^\omega$, $Y\setminus X$ is finite, then $Y\setminus A$ is finite. This shows that $\mathcal{A}_D$ is open. To show density, let $X\in[\omega]^\omega$ be a set not in $\mathcal{A}_D$ (the case where $X\in\mathcal{A}_D$ is trivial.) Since $\mathcal{A}$ is MAD, there is a $Y\in\mathcal{A}$ such that $X\cap Y$ is infinite. Then $X\cap Y\in\mathcal{A}_D$ and $X\cap Y\subseteq X$, proving density. 

If $y\in \mathcal{Y}$ is any arbitrary element, we want to show that $y\not\in\mathcal{A_D}$. For contradiction suppose that $y\in\mathcal{A}_D$, then there is some $A\in\mathcal{A}$ such that $y\setminus A$ is finite. By being MAD, if $y\in\mathcal{A}_D$ then for any distinct $A_0,A_1\in\mathcal{A}$ either $A_0\cap y$ or $A_1\cap y$ is finite. By Lemma \ref{MAD lemma} this implies there exists some $A\in\mathcal{A}$ such that $y\setminus A$ is finite, which means that $y\in\mathcal{A}_D$. By contradiction, no element of $\mathcal{Y}$ is an element of $\mathcal{A}_D$.
\end{proof}
\end{theorem}

In light of Theorem \ref{contiuumn theorem}, Theorem \ref{max(r,d)=hom}, and Proposition \ref{r=hom par=s}, the dual realization of the theorem below proves that most of the cardinal characteristics discussed in the thesis have cardinality less than or equal to $\mathfrak{c}$. The ones not included are those not mapped to by $\mathfrak{d}$. In addition, the theorem below shows that $\mathfrak{h}$ functions as a lower bound for all the cardinals discussed in this section.

\begin{theorem}$\mathfrak{h}\leq\mathfrak{par_2}$\label{h<par}
\begin{proof}
We construct a morphism $\varphi:([\omega]^{\omega},P_2,\neg H^*)\rightarrow([\omega]^{\omega},DO,\not\in)$

Let $\varphi_{-}:[\omega]^{\omega}\rightarrow[\omega]^{\omega}$ be the identity function. For any 2-coloring $\pi$ let $\varphi_{+}(\pi)$ be the family $F_{\pi}\subseteq[\omega]^{\omega}$  of all infinite sets which are almost homogeneous for $\pi$. If $N\in[\omega]^{\omega}$ is an arbitrary infinite set that is not almost homogeneous $\pi$, then $N\not\in F_{\pi}$.

We must now verify that $F_{\pi}$ is a dense-open subset. Density follows directly from Ramsey's theorem. To verify that $F_{\pi}$ is open, observe that adding finitely many elements to an almost homogeneous set results in another almost homogeneous set.
\end{proof}  
\end{theorem}

\begin{theorem}
$\mathfrak{h}\leq\mathfrak{g}\leq\mathfrak{d}$
\begin{proof}
We can construct morphism $\varphi_1:([\omega]^\omega,G,\notin)\rightarrow([\omega]^\omega,\mathrm{DO},\notin)$,  by letting both $\varphi_{1_+}$ and $\varphi_{1_-}$ be the identity. The only non-trivial aspect is showing that any groupwise dense set is dense-open. By definition, any groupwise dense set is open. To show density, let $\mathcal{G}$ be an arbitrary groupwise dense family. For any $X\in[\omega]^\omega$, let $\Pi_X\in\mathrm{IP}$ be such that each interval contains at least one element from $X$. By closure under almost subsets, the union over each $I\in\Pi_X$ gives an infinite subset of $X$ which is in $\mathcal{G}$.

To produce a morphism $\varphi_2:(\prescript{\omega}{\hphantom{.}}{\omega},\prescript{\omega}{\hphantom{.}}{\omega},<^*)\rightarrow([\omega]^\omega,G,\notin)$, we let $\varphi_{2_-}$ associate $A\in[\omega]^\omega$ with the function $A_f\in \prescript{\omega}{\hphantom{.}}{\omega}$ enumerating it. Assume that $A_f<^*g$ and let $\varphi_{2_+}(g)=F_g$, defined as: $$F_g:=\{x\in[\omega]^\omega:\text{ there exists infintely many }n\in\omega\text{ such that }[n,f(n))=\emptyset\}.$$

If, for contradiction, we supposed that $A\in F_g$, this would imply that there are infinitely many $n\in\omega$ such that $A\cap[n,f(n))=\emptyset$. If $n_0\in\omega$ is the point past which $g$ dominates $f$, then for each $n>n_0$, $n\leq a_k\leq g(k)$. This implies that for all but finitely many $n\in\omega$, $A\cap[n,f(n))\neq\emptyset$.
\end{proof}
\end{theorem}
\begin{figure}[ht]
\centering
\[
\xymatrix@R=1.02em@C=2em{
   & \mathfrak{c} \ar[d] & \\
   & \mathfrak{hom}_2 \ar[dl] \ar[d] & \\
 \mathfrak{r} \ar[d] & \mathfrak{d} \ar[d] \ar[dr]\ar[dl] & \\
 \mathfrak{b} \ar[dr] & \mathfrak{s} \ar[d] & \mathfrak{g} \ar[ddl]\\
 & \mathfrak{par}_2 \ar[d] & \\
 & \mathfrak{h} \ar[d] & \\
 & \omega_1 &
}
\]
\caption{A Diagram of the Additional Cardinal Characteristics}
\label{fig:extra_diagram}
\end{figure}
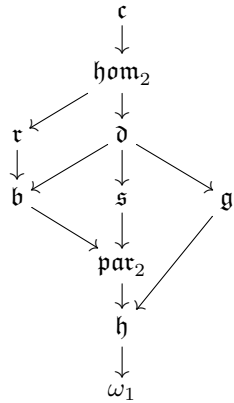

\end{document}